\documentclass[notitlepage,leqno,10pt]{article}

\usepackage{latexsym}
\usepackage{amsmath}
\usepackage{amsfonts}
\usepackage{amssymb}
\renewcommand{\a }{\alpha }
\renewcommand{\b }{\beta }
\renewcommand{\d}{\delta }
\newcommand{\D }{\Delta }

\newcommand{\e }{\varepsilon }
\newcommand{\g }{\gamma}

\renewcommand{\l }{\lambda }

\newcommand{\n }{\nabla }
\newcommand{\var }{\varphi }

\newcommand{\s }{\sigma }
\newcommand{\Sig }{\Sigma}

\renewcommand{\th }{\theta }
\renewcommand{\o }{\omega }

\newcommand{\ov}{\overline}
\newcommand{\intbar}{\mathop{\int\makebox(-13.5,0){\rule[4pt]{.7em}{0.3pt}}%
\kern-6pt}\nolimits}

\newcommand{\be}{\begin{equation}}
\newcommand{\ee}{\end{equation}}
\newenvironment{pf}{\noindent{\sc Proof}.\enspace}{\rule{2mm}{2mm}\medskip}
\newenvironment{pfn}{\noindent{\sc Proof}}{\rule{2mm}{2mm}\medskip}

\newcommand{\R}{\mathbb{R}}

\newcommand{\N}{\mathbb{N}}

\author{Andrea MALCHIODI}

\date{}

\title{Compactness of solutions to some geometric \\
fourth-order equations}

\begin{document}

\newtheorem{lem}{Lemma}[section]
\newtheorem{pro}[lem]{Proposition}
\newtheorem{thm}[lem]{Theorem}
\newtheorem{rem}[lem]{Remark}
\newtheorem{cor}[lem]{Corollary}
\newtheorem{df}[lem]{Definition}

\maketitle

\begin{center}

{\small Sissa, Via Beirut 2-4, 34014 Trieste, Italy}

\end{center}

\footnotetext[1]{E-mail address: malchiod@sissa.it}

\

\

\noindent {\sc abstract}. We prove compactness of solutions to
some fourth order equations with exponential nonlinearities on
four manifolds. The proof is based on a refined bubbling analysis,
for which the main estimates are given in integral form. Our
result is used in a subsequent paper to find critical points (via
minimax arguments) of some geometric functional, which give rise
to conformal metrics of constant $Q$-curvature. As a byproduct of
our method, we also obtain compactness of such metrics.

\begin{center}

\bigskip\bigskip

\noindent{\it Key Words:} Fourth order equations, Blow-up
analysis, Geometric PDEs

\bigskip

\centerline{\bf AMS subject classification: 35B33, 35J35, 53A30,
53C21}

\end{center}

\section{Introduction}

\noindent Consider a compact four-dimensional manifold $(M, g)$
with Ricci tensor $Ric_g$ and scalar curvature $R_g$. The $Q$-{\em
curvature} and the {\em Paneitz operator}, introduced in
\cite{bo}, \cite{p1} and \cite{p2}, are defined respectively by
\begin{equation}\label{eq:Q}
    Q_g = - \frac{1}{12} \left( \D_g R_g - R_g^2 + 3 |Ric_g|^2
    \right);
\end{equation}
\begin{equation}\label{eq:P}
    P_g(\var) = \D_g^2 \var + \hbox{div} \left( \frac 23 R_g g - 2 Ric_g
    \right) d \var,
\end{equation}
where $\var$ is any smooth function on $M$, see also the survey
\cite{cy99}.

The $Q$-curvature and the Paneitz operator arise in several
contexts in the study of four-manifolds and of particular interest
is their role, and their mutual relation, in conformal geometry.
In fact, given a metric $\tilde{g} = e^{2 w} g$, the following
equations hold
\begin{equation}\label{eq:confP}
    P_{\tilde{g}} = e^{- 4 w} P_g; \qquad \qquad \qquad
    P_g w + 2 Q_g = 2 Q_{\tilde{g}} e^{4 w}.
\end{equation}
A first connection to the topology of a manifold is a Gauss-Bonnet
type formula. If $W_g$ denotes the Weyl's tensor of $M$, then one
has
$$
    \int_M \left( Q_g + \frac{|W_g|^2}{8} \right) dV_g = 4 \pi^2
    \chi(M),
$$
where $dV_g$ stands for the volume element in $(M,g)$ and
$\chi(M)$ is the Euler characteristic of $M$. In particular, since
$|W_g|^2$ is a pointwise conformal invariant, it follows that
$\int_M Q_g dV_g$ is a {\em global} conformal invariant.

To mention some geometric applications we recall three results
proven by Gursky, \cite{g2}, and by Chang, Gursky and Yang,
\cite{cgyann}, \cite{cgy} (see also \cite{g}). If a manifold of
positive Yamabe class satisfies $\int_M Q_g dV_g > 0$, then its
first Betti number vanishes. Moreover there exists a conformal
metric with positive Ricci tensor, and hence $M$ has finite
fundamental group. Furthermore, under the additional quantitative
assumption $\int_M Q_g dV_g > \frac 18 \int_M |W_g|^2 dV_g$, $M$
bust be diffeomorphic to the four-sphere or to the projective
space. In particular the last result is a conformally invariant
improvement of a theorem by Margerin, \cite{marrg}, which assumed
pointwise pinching conditions on the Ricci tensor in terms of
$W_g$.

Finally, we also point out that the Paneitz operator and the
$Q$-curvature (together with their higher-dimensional analogues,
see \cite{br}, \cite{br2}, \cite{fg}, \cite{gjms}) appear in the
study of Moser-Trudinger type inequalities, log-determinant
formulas and the compactification of locally conformally flat
manifolds, see \cite{beck}, \cite{bo}, \cite{bcy}, \cite{cqy},
\cite{cqy2}, \cite{cy95}.

\

\noindent As for the uniformization theorem, one can ask whether
every four-manifold $(M,g)$ carries a conformal metric $\tilde{g}$
for which the corresponding $Q$-curvature $Q_{\tilde{g}}$ is a
constant. Writing $\tilde{g} = e^{2 w} g$, by \eqref{eq:confP} the
problem is equivalent to finding a solution of the equation
\begin{equation}\label{eq:Qc}
    P_g w + 2 Q_g = 2 \ov{Q} e^{4 w},
\end{equation}
where $\ov{Q}$ is a real constant. In view of the regularity
results in \cite{uv}, solutions of \eqref{eq:Qc} can be found as
critical points of the following functional
\begin{equation}\label{eq:II}
    II(u) = \langle P_g u, u \rangle + 4 \int_M Q_g u dV_g - k_P
    \log \int_M e^{4u} dV_g; \qquad u \in H^{2}(M),
\end{equation}
where we are using the notation
$$
  \langle P_g u, v \rangle = \int_M \left( \D_g u \D_g v + \frac 23 R_g \n_g u
  \cdot \n_g v - 2 (Ric_g \n_g u, \n_g v) \right) dV_g; \qquad u, v \in
  H^2(M),
$$
and where
\begin{equation}\label{eq:kp}
    k_P = \int_M Q_g dV_g.
\end{equation}
Problem \eqref{eq:Qc} has been solved in \cite{cy95} for the case
in which $P_g$ is a positive operator and $k_P < 8 \pi^2$ ($8
\pi^2$ is the value of $k_P$ on the standard sphere). Under these
assumptions by the Adams inequality, see \eqref{eq:ada2}, the
functional $II$ is bounded from below and coercive, hence
solutions can be found as global minima. The result has also been
extended in \cite{b1} to higher-dimensional manifolds (regarding
higher-order operators and curvatures) using a geometric flow. A
first sufficient condition to ensure these hypotheses was given by
Gursky in \cite{g2}. He proved that if the Yamabe invariant is
positive and if $k_P > 0$, then $P_g$ is positive definite and
moreover $k_P \leq 8 \pi^2$, with the equality holding if and only
if $M$ is conformally equivalent to $S^4$. Other more general
sufficient conditions are given in \cite{gv}. The solvability of
\eqref{eq:Qc} also turns out to be useful in the study of some
interesting class of fully non-linear equations, as it has been
shown in \cite{cgy}, with the remarkable geometric consequences
mentioned above.

\

\noindent We are interested here in the more general case when
$P_g$ has no kernel and $k_P \neq 8 k \pi^2$ for $k = 1, 2,
\dots$. These conditions are generic, and in particular include
manifolds with negative curvature or negative Yamabe class, for
which $k_P$ can be bigger than $8 \pi^2$.

In the case under investigation the functional $II$ can be
unbounded from below, and hence it is necessary to find extrema
which are possibly saddle points. As we shall explain below, in
order to find these critical points it is useful to study
compactness of solutions to perturbations of \eqref{eq:Qc}.

Therefore we consider the following sequence of problems
\begin{equation}\label{eq:pl}
    P_g u_l + 2 Q_l = 2 k_l e^{4 u_l} \qquad \hbox{ in } M,
\end{equation}
where $(k_l)_l$ are constants and where
\begin{equation}\label{eq:Ql}
    Q_l \to Q_0 \qquad \hbox{ in } C^0(M).
\end{equation}
Without loss of generality, we can assume that the sequence
$(u_l)_l$ satisfies the volume normalization
\begin{equation}\label{eq:noul}
    \int_M e^{4 u_l} dV_g = 1 \qquad \qquad \hbox{ for all } l,
\end{equation}
which implies that we must choose $k_l = \int_M Q_l dV_g$.

\

\noindent Our main result is the following.

\begin{thm}\label{th:bd}
Suppose $ker \; P_g = \{constants\}$ and that $(u_l)_l$ is a
sequence of solutions of \eqref{eq:pl}, \eqref{eq:noul}, with
$(Q_l)_l$ satisfying \eqref{eq:Ql}. Assume also that
\begin{equation}\label{eq:kp2}
    k_0 := \int_M Q_0 dV_g \neq 8 k \pi^2 \qquad \qquad
    \hbox{ for } k = 1, 2, \dots.
\end{equation}
Then $(u_l)_l$ is bounded in $C^\a(M)$ for any $\a \in (0,1)$.
\end{thm}

\noindent The main application of Theorem \ref{th:bd} concerns the
case $Q_0 = Q_g$. Indeed, if a sequence of solutions to
\eqref{eq:pl}-\eqref{eq:noul} can be produced, its weak limit will
be a critical point of the functional $II$ and a solution of
\eqref{eq:Qc}. This is indeed verified in \cite{dm2} under the
assumptions of Theorem \ref{th:bd} (with $Q_0 = Q_g$). As a
consequence one finds conformal metrics with constant
$Q$-curvature for a large class of four manifolds. We have indeed
the following result, announced in the preliminary note \cite{dm1}
with some sketch of the ideas.

\begin{thm} (\cite{dm2})
Suppose $ker \; P_g = \{constants\}$, and assume that $k_P \neq 8
k \pi^2$ for $k = 1, 2, \dots$. Then equation \eqref{eq:Qc} has a
solution.
\end{thm}

\noindent The proof requires a minimax scheme which becomes more
and more involved as $k$ increases and when $P_g$ possesses
negative eigenvalues. This scheme extends that in \cite{djlw},
which in our case would correspond to $P_g \geq 0$ and $k_0 \in (8
\pi^2, 16 \pi^2)$.

The way we use Theorem \ref{th:bd} in \cite{dm2} is the following.
First, for $\rho$ in a neighborhood of $1$, we introduce the
modified functional
$$
    II_\rho(u) = \langle P_g u, u \rangle + 4 \rho \int_M Q_g u
    dV_g - k_P \rho \log \int_M e^{4u} dV_g; \qquad u \in H^{2}(M),
$$
and, using the minimax scheme, we prove existence of Palais-Smale
sequences at some level $c_\rho$. It turns out that the function
$\rho \mapsto c_\rho$ is a.e. differentiable and, following an
idea in \cite{str} (used also in \cite{djlw}, \cite{jt},
\cite{st}), we prove existence of critical points of $II_\rho$ for
those values of $\rho$ at which $c_\rho$ is differentiable. Then
we are led to consider \eqref{eq:pl} for $Q_l = \rho_l Q_g$, where
$(\rho_l)_l$ is a suitable sequence.

\

\noindent Theorem \ref{th:bd} applies also to any sequence of
smooth solutions of \eqref{eq:Qc}. Therefore, as another
application, we have the following result, which extends a
compactness theorem in \cite{cy95}.

\begin{cor}\label{th:comp}
Suppose $ker \; P_g = \{constants\}$ and that $k_p \neq 8 k \pi^2$
for $k = 1, 2, \dots$. Suppose $(u_l)_l$ is a sequence of
solutions of \eqref{eq:Qc} satisfying \eqref{eq:noul}. Then, for
any $m \in \N$, $(u_l)_l$ is bounded in $C^m(M)$.
\end{cor}

\noindent Corollary \ref{th:comp} has a counterpart in \cite{li}
(see also \cite{cl2}), where compactness of solutions is proved
for a mean field equation on compact surfaces.

\

\noindent The case when $k_P$ is an integer multiple of $8 \pi^2$
is more delicate, and should require an asymptotic analysis as in
\cite{bah}, \cite{cl1}, \cite{cl2}, \cite{li} (see also the
references therein). An interesting particular case of this
situation is the standard sphere. Being an homogeneous space, the
$Q$-curvature is already constant and indeed all the solutions of
\eqref{eq:Qc} on $S^4$, which have been classified in \cite{cy97},
arise from conformal factors of M\"obius transformations.
Henceforth, a natural problem to consider is to prescribe the
$Q$-curvature as a given function $f$ on $S^4$. Some results in
this direction are given in \cite{b2}, \cite{ms} and \cite{wx}.
Typically, the methods are based on blow-up or asymptotic analysis
combined with Morse theory, in order to deal with a possible loss
of compactness.

\

\noindent The Paneitz operator and the $Q$-curvature can be
considered as natural extensions to four-manifolds of,
respectively, the Laplace Beltrami operator $\D_g$ and the Gauss
curvature $K_g$ on two-dimensional surfaces. In fact, similarly ro
$P_g$ and $Q_g$, these transform according to the equations
\begin{equation}\label{eq:ctl}
    \D_{\tilde{g}} = e^{-2 w} \D_g; \qquad \qquad
    - \D_g w + K_g = K_{\tilde{g}} e^{2w},
\end{equation}
where, again, $\tilde{g} = e^{2 w} g$. Hence, in the case of a
flat domain $\Omega \subseteq \R^2$, one is led to study equations
of the form
\begin{equation}\label{eq:K(x)}
    - \D v_l = K_l(x) e^{2 v_l} \qquad \hbox{ in } \Omega.
\end{equation}

In \cite{bm} the authors proved, among other things, that if
$(K_l)_l$ are non-negative, uniformly bounded in
$L^\infty(\Omega)$ and if $\int_\Omega e^{2 u_l} \leq C$, then
either $(v_l)_l$ stays bounded in $L^\infty_{loc}(\Omega)$, or
$v_l \to - \infty$ on the compact subsets of $\Omega$, or $K_l
e^{2 v_l}$ concentrates at a finite number of points in $\Omega$,
namely $K_l e^{2 v_l} \rightharpoonup \sum_{i=1}^j \a_i \d_{x_i}$
($\d_{x_i}$ stands for the Dirac mass at $x_i$). In the latter
case, they also proved that each $\a_i$ is greater or equal than
$4 \pi$. This result was specialized in \cite{ls} where, assuming
that $K_l \to K_0$ in $C^0(\ov{\Omega})$ and using the sup+inf
inequalities in \cite{bls}, \cite{sh}, the authors proved that
each $\a_i$ is indeed an integer multiple of $4 \pi$. Chen showed
then in \cite{chen} that the case of a multiple bigger than $1$
may indeed occur. On the other hand, if $\Omega$ is replaced by a
compact surface (subtracting a constant term to the right-hand
side, to get solvability of the equation), then each $\a_i$ is
precisely $4 \pi$, see \cite{li}. The same result is obtained in
\cite{os} for approximate solutions in domains, but with an extra
assumption on the $L^\infty$ norm of the error terms.

\

\noindent Our argument for the proof of Theorem \ref{th:bd}, which
we outline below, relies on proving a quantization result for the
volume of blowing-up solutions as in \cite{ls}. The main idea is
to show that at every blow-up point the volume is a multiple of
$\frac{8 \pi^2}{k_0}$. Then, proving also that there is no
residual volume amount, we reach a contradiction with
\eqref{eq:noul} since we are assuming that $k_0$ is not an integer
multiple of $8 \pi^2$. However, instead of using pointwise
estimates on the solutions, as in \cite{bm} or \cite{ls}, our
results are mainly given in integral form, see Remark \ref{r:ps}.

Except for the last subsection, we work under the assumption
\begin{equation}\label{eq:kp3}
    k_0 \in (8 k \pi^2, 8 (k+1) \pi^2),
    \quad k \in \N,
\end{equation}
since this case contains most of the difficulties.

The plan of the paper (and the strategy of the proof) is the
following. In Section \ref{s:pr} we collect some preliminary facts
including a modified version of the Adams inequality, to deal with
the presence of negative eigenvalues, and some $L^p$ estimates on
the first, second and third derivatives of the solutions.

In Section \ref{s:bub} we derive a compactness criterion based on
the amount of concentration of the nonlinear term, see Proposition
\ref{p:cc}, and then we study the asymptotic profile of $u_l$ near
the concentration points. In particular we prove that the minimal
volume accumulation is $\frac{8 \pi^2}{k_0}$, see
\eqref{eq:intalrlu}.

In Section \ref{s:sbu}, which is the core of our analysis, we
introduce the notion of {\em simple blow-up} (adopting the
terminology used by R.Schoen) and we show in Proposition
\ref{p:sbu} that at such blow-ups the accumulation is exactly
$\frac{8 \pi^2}{k_0}$. In order to prove this we use some integral
form of the Harnack inequality, see in particular Subsection
\ref{ss:ih}, combined with a careful ODE analysis for the function
$r \mapsto \ov{u}_{r,l}$. Here $\ov{u}_{r,l}$ denotes, naively,
the average of $u_l$ on an annulus $A_r$ of radii $r$ and $2 r$
centered near a concentration point.

Finally, in Section \ref{s:fa} we show how a general blow-up
situation can be essentially reduced to the case of finitely-many
simple blow-ups. In particular, we prove that at any general
blow-up point the amount of concentration is an integer multiple
of $\frac{8 \pi^2}{k_0}$. Recalling the normalization
\eqref{eq:noul} and that $k_0 \neq 8 k \pi^2$ for any integer $k$,
we reach then a contradiction to the fact that $(u_l)_l$ is
unbounded in some $C^\a$ norm. In Subsection \ref{ss:pnot>0} we
consider the case $k_0 < 8 \pi^2$, which is easier and requires
only the analysis of Section \ref{s:bub}.

In our proof we exploit crucially the fact that we are working on
a compact manifold, since we often make use of the Green's
representation formula. We also point out that our assumptions on
$M$ are generic and do not require the metric to be locally
conformally flat or Einstein.

\begin{rem}\label{r:ps}
It is an open problem to understand whether the functional $II$
itself (see \eqref{eq:II}) possesses bounded Palais-Smale
sequences, or equivalently if it is possible to find solutions of
\eqref{eq:Qc} without introducing the perturbed functional
$II_\rho$.

The reason why we kept most of our estimates in integral form is
that many of them could be applied to functions of class $H^2$
only (not necessarily smooth or bounded) and we hope that some
could be useful to understand the question. At the moment, in
particular, the counterparts of Proposition \ref{p:sbu} is missing
for Palais-Smale sequences and we need the full rigidity of
equation \eqref{eq:pl}. For related topics see \cite{os}.
\end{rem}

\

\begin{center}

{\bf Acknowledgements}

\end{center}

\noindent This work was started when the author was visiting IAS
in Princeton, and continued during his stay at ETH in Z\"urich,
Laboratoire Jacques-Louis Lions in Paris, Sissa in Trieste and IMS
in Singapore. He is very grateful to all these institutions for
their kind hospitality. The author has been supported by
M.U.R.S.T.  under the national project {\em Variational methods
and nonlinear differential equations}, and by the European Grant
ERB FMRX CT98 0201.

\section{Notation and preliminaries}\label{s:pr}

In this brief section we collect some useful preliminary facts,
and in particular we state a version of the Moser-Trudinger
inequality involving the Panetiz operator. In the following
$B_r(p)$ stands for the metric ball of radius $r$ and center $p$.
We also denote by $|x - y|$ the distance of two points $x, y \in
M$. $H^2(M)$ is the Sobolev space of functions on $M$ which are in
$L^2(M)$ together with their first and second derivatives. Large
positive constants are always denoted by $C$, and the value of $C$
is allowed to very from formula to formula and also within the
same line.

\

\noindent As already mentioned, throughout most of the paper we
will work under the assumption \eqref{eq:kp3}. When the operator
$P_g$ is positive definite, by the Poincar\'e inequality the $H^2$
norm is equivalent to the following one
\begin{equation}\label{eq:nuP}
    \|u\|^2 = \langle P_g u, u \rangle + \int_M u^2 dV_g; \qquad u \in
    H^2(M).
\end{equation}
Being $M$ four-dimensional, $H^2(M) \hookrightarrow L^p(M)$ for
all $p > 1$. We have indeed the following limit-case embedding,
proved in \cite{ada} and \cite{bcy} for the operator $\D^2$ and
extended in \cite{cy95} for the Paneitz operator.

\begin{pro}\label{p:ada}
If $P_g \geq 0$, there exists a positive constant $C$ depending on
$M$ such that
\begin{equation}\label{eq:ada1}
    \int_M e^{\frac{32 \pi^2 (u - \ov{u})^2}{\langle P_g u, u \rangle}}
    dV_g \leq C; \qquad \hbox{ for every } u \in H^{2}(M),
\end{equation}
where $\ov{u} = \frac{1}{Vol(M)} \int_M u dV_g$ denotes the
average of $u$ on $M$. The last formula implies
\begin{equation}\label{eq:ada2}
    \log \int_M e^{4(u - \ov{u})} dV_g \leq C + \frac{1}{8 \pi^2}
    \langle P_g u, u \rangle; \qquad \hbox{ for every } u \in
    H^{2}(M).
\end{equation}
\end{pro}

\

\noindent Here we are interested in the case in which $P_g$ might
possess some negative eigenvalues. We denote by $V \subseteq
H^2(M)$ the direct sum of the eigenspaces corresponding to
negative eigenvalues of $P_g$. Of course the dimension of $V$ is
finite, say $\ov{k}$, and since $P_g$ has no kernel and is
self-adjoint we can find an orthonormal basis of eigenfunctions
$\hat{v}_1, \dots, \hat{v}_{\ov{k}}$ of $V$ with the properties
\begin{equation}\label{eq:hatv1k}
    P_g \hat{v}_i = \l_i \hat{v}_i, \quad i = 1, \dots, \ov{k};
  \qquad \quad \l_1 \leq \l_2 \leq \dots \leq \l_{\ov{k}} < 0 <
  \l_{\ov{k}+1} \leq \dots,
\end{equation}
where the $\l_i$'s are the eigenvalues of $P_g$. Having introduced
the subspace $V$, we need a modified version of the Adams
inequality.

\begin{lem}\label{l:adaneg}
Suppose $P_g$ possesses some negative eigenvalues, that $ker P_g =
\{constants\}$, and let $V$ denote the direct sum of the negative
eigenspaces of $P_g$. Then there exists a constant $C$ such that
\begin{equation}\label{eq:adaneg0}
    \int_M e^{\frac{32 \pi^2 (u - \ov{u})^2}{\langle P_g u, u \rangle}}
    dV_g \leq C; \qquad \hbox{ for every function } u \in H^2(M)
    \hbox{ with } \hat{u} = 0.
\end{equation}
Here $\hat{u}$ denotes the component of $u$ in $V$. As a
consequence one has
\begin{equation}\label{eq:adaneg}
  \log \int_M e^{4(u - \ov{u})} dV_g \leq C + \frac{1}{8 \pi^2}
    \langle P_g u, u \rangle, ; \qquad \hbox{ for every function } u
    \in H^2(M) \hbox{ with } \hat{u} = 0.
\end{equation}
\end{lem}

\begin{pf}
The proof is a variant of the arguments of \cite{bcy} and
\cite{cy95}. If $\hat{v}_1, \dots, \hat{v}_{\ov{k}}$ and $\l_1,
\dots, \l_{\ov{k}}$ are as in \eqref{eq:hatv1k}, we introduce the
following positive-definite pseudo-differential operator $P^+_g$
$$
  P^+_g u = P_g u - 2 \sum_{i=1}^{\ov{k}} \l_i \left( \int_M u \hat{v}_i
  dV_g \right) \hat{v}_i.
$$
Basically, we are reversing the sign of the negative eigenvalues
of $P_g$. The operator $P_g^+$ admits the following Green's
function
$$
  G^+(x,y) = G(x,y) - 2 \sum _{i=1}^{\ov{k}} \l_i \hat{v}_i(x)
  \hat{v}_i(y),
$$
where $G(x,y)$ corresponds to $P_g$. Then the arguments of
\cite{cy95} (see also \cite{ada}, \cite{bcy}), which are based on
representations for pseudo-differential operators, can be adapted
to the case of $P_g^+$, yielding
$$
\int_M e^{\frac{32 \pi^2 (u - \ov{u})^2}{\langle P^+_g u, u
\rangle}} dV_g \leq C \qquad \qquad \hbox{ for every } u \in
H^2(M).
$$
Applying the last formula to functions for which $\hat{u} = 0$, we
obtain \eqref{eq:adaneg0}. Finally, from the elementary inequality
$4 a b \leq 32 \pi^2 a^2 + \frac{1}{8 \pi^2} b^2$, applied with $a
= (u - \ov{u})$ and $b = \langle P_g u, u \rangle$, we also obtain
\eqref{eq:adaneg}.
\end{pf}

\

\noindent Theorem \ref{th:bd} is proved by contradiction. We claim
that unboundedness in some $C^\a$ norm is equivalent (under the
assumption \eqref{eq:kp3}, which implies $k_l > 0$ for $l$ large)
to the following condition
\begin{equation}\label{eq:unbdul}
    \|u_l - \ov{u}_l\| \to + \infty \qquad \hbox{ as } l
    \to + \infty.
\end{equation}
In order to prove this we first notice that, by \eqref{eq:noul}
and the Jensen inequality, $\ov{u}_l$ is uniformly bounded from
above. Assuming that $\|u_l - \ov{u}_l\|$ is uniformly bounded
(which implies, in the above notation, that also $\|u_l - \ov{u}_l
- \hat{u}_l\|$ is uniformly bounded), then by \eqref{eq:adaneg}
the right-hand side of \eqref{eq:pl} is also uniformly bounded in
$L^p(M)$ for every $p > 1$. By elliptic regularity, then $(u_l)_l$
would be uniformly bounded in $W^{4,p}(M)$, and hence in $C^\a(M)$
for any $\a \in (0,1)$ by the Sobolev embeddings.

Hence from now on we assume that there exists a sequence $(u_l)_l$
satisfying \eqref{eq:pl}-\eqref{eq:noul} and \eqref{eq:unbdul}.

\

\noindent We prove now a preliminary integrability result on the
first, second and third derivatives of $u_l$.

\begin{lem}\label{l:lodlp}
Let $(u_l)_l$ be a sequence of solutions of
\eqref{eq:pl}-\eqref{eq:noul}, with $(Q_l)_l$ satisfying
\eqref{eq:Ql}, and let $p \geq 1$. Then there is a constant $C$
depending only on $p$, $M$ and $k_0$ such that, for $r$
sufficiently small and for any $x \in M$ there holds
$$
  \int_{B_r(x)} |\n^3 u_l|^p dV_g \leq C r^{4 - 3p}; \qquad \int_{B_r(x)}
  |\n^2 u_l|^p dV_g \leq C r^{4 - 2p}; \qquad \int_{B_r(x)} |\n
  u_l|^p dV_g \leq C r^{4 - p},
$$
where, respectively, $p < \frac 43, p < 2$ and $p < 4$.
\end{lem}

\begin{pf}
We write
$$
  P_g u_l = f_l,
$$
where
\begin{equation}\label{eq:fl}
    f_l = 2 k_l e^{4 u_l} - 2 Q_l.
\end{equation}
We have the following representation formula
\begin{equation}\label{eq:Grepwl}
    u_l(x) = \ov{u}_l + \int_M G(x,y) f_l(y)
    dV_g(y), \qquad \hbox{ for a.e. } x \in M,
\end{equation}
where, by the results in \cite{cy95}, $G : M \times M \setminus
\hbox{diag}$ is symmetric and satisfies
\begin{equation}\label{eq:G}
    \left| G(x,y) - \frac{1}{8 \pi^2} \log \frac{1}{|x - y|}
    \right| \leq C, \qquad \qquad x, y \in M, x \neq y,
\end{equation}
while for its derivatives there holds
\begin{equation}\label{eq:derG}
  |\n G(x,y)| \leq C \frac{1}{|x-y|}; \qquad |\n^2 G(x,y)| \leq C
  \frac{1}{|x-y|^2}; \qquad |\n^3 G(x,y)| \leq C \frac{1}{|x-y|^3}.
\end{equation}
The last two estimates in \eqref{eq:derG} are not shown in
\cite{cy95} but they can be derived with the same approach, by an
expansion of $G$ at higher order using the parametrix, see also
\cite{aul}. Similarly (this formula will be used later in the
paper), one also finds that
\begin{equation}\label{eq:derG2}
    \n_x G(x,y) = \frac{1}{8 \pi^2} \n_x \log \frac{1}{|x - y|} + O(1).
\end{equation}

Recalling the definition of $f_l$ in \eqref{eq:fl}, we obtain
$$
  |\n^3 u_l|(x) \leq C \int_M \frac{1}{|x-y|^3} |f_l(y)| dV_g(y), \qquad
  \hbox{ for a.e. } x \in M.
$$
Then, from the Jensen's inequality it follows that
$$
  |\n^3 u_l|^p(x) \leq C \int_M \left( \frac{\|f_l\|_{L^1(M)}}{|x-y|^3}
  \right)^p \frac{|f_l(y)|}{\|f_l\|_{L^1(M)}} dV_g(y) \qquad
  \hbox{ for a.e. } x \in M.
$$
The Fubini's Theorem implies
$$
  \int_{B_x(\ov{x})} |\n^3 u_l|^p(x) dV_g(x) \leq C \sup_{y \in M}
  \int_{B_x(\ov{x})} \frac{1}{|x-y|^{3p}} dV_g(x) \leq C \int_{B_x(\ov{x})}
  \frac{1}{|x-\ov{x}|^{3p}} dV_g(x).
$$
The last integral is finite provided $3 p < 4$, as in our
assumptions, and can be estimated using polar coordinates, giving
$$
  \int_{B_x(\ov{x})} |\n^3 u_l|^p(x) dV_g(x) \leq C(p,M) r^{4-3p}.
$$
This concludes the proof of the first inequality in the statement
of the lemma. The remaining two follow similarly.
\end{pf}

\section{The bubbling phenomenon}\label{s:bub}

In this section we study the local behavior of unbounded sequences
of solutions at a concentration point. In subsection \ref{ss:cc}
we give compactness criteria when the amount of concentration is
below a certain threshold. Then, in Subsection \ref{ss:ap}, we
reduce ourselves to the preceding situation using a scaling
argument. As a byproduct we describe the asymptotic profile of
$u_l$, proving that it has the form of a {\em standard bubble},
and we show that the amount of volume concentration at any blow-up
point is greater or equal than $\frac{8 \pi^2}{k_0}$.

\subsection{Concentration-compactness}\label{ss:cc}

In this subsection we give a concentration-compactness criterion
for solutions of the equation $P_g v = h$ on $M$. In the case of
the sphere a similar result has been shown in \cite{b1}, and our
proof basically goes along the same line. However we prefer to
write the details, since some of them will be needed in the
following.

\begin{pro}\label{p:cc}
Let $(h_l)_l \subseteq L^1(M)$ be a sequence of functions
satisfying $\int_M |h_l| dV_g \leq C$ for every $l$. Let $v_l$ be
solutions of $P_g v_l = h_l$ on $M$. Then, up to a subsequence,
either for every $l$
$$
  \int_M e^{\a (v_l - \ov{v}_l)} dV_g \leq C \qquad \hbox{ for some }
  C > 0 \hbox{ and some } \a > 4,
$$
or there exists points $x_1, \dots, x_L \in M$ such that, for any
$r > 0$ and any $i \in \{1, \dots, L\}$ there holds
\begin{equation}\label{eq:conc}
    \liminf_{l \to + \infty} \int_{B_r(x_i)} |h_l| dV_g \geq 8 \pi^2.
\end{equation}
\end{pro}

\begin{pf}
Assume the second alternative does not occur, namely
\begin{equation}\label{eq:concontr}
    \hbox{ for every } x \in M \hbox{ there exists } r_x > 0
    \hbox{ such that } \int_{B_{r_x}(x)} |h_l| dV_g \leq
    8 \pi^2 - \d_x,
\end{equation}
for some $\d_x > 0$ and for $l$ sufficiently large. We cover $M$
with $j$ balls $B_i := B_{\frac{r_{x_i}}{2}}(x_i)$, $i = 1, \dots,
j$. Using \eqref{eq:Grepwl} and setting $B_{r_{x_i}}(x_i) =
\tilde{B}_i$, for a.e. $x \in B_i$ we can write
\begin{equation}\label{eq:usum}
    v_l(x) - \ov{v}_l = \int_{\tilde{B}_i} h_l(y) G(x,y) dV_g(y)
    + \int_{M \setminus \tilde{B}_i} h_l(y) G(x,y) dV_g(y).
\end{equation}
Hence if $\a > 0$, for a.e. $x \in B_i$ we have
\begin{equation}\label{eq:expau}
    \exp \left[\a (v_l(x) - \ov{v}_l) \right] = \exp \left[ \int_{\tilde{B}_i} \a G(x,y) h_l(y)
   dV_g(y) \right] \exp \left[ \int_{M \setminus \tilde{B}_i}
   \a G(x,y) h_l(y) dV_g(y) \right].
\end{equation}
Since $G$ is smooth outside the diagonal, and since $\int_M |h_l|
dV_g$ is uniformly bounded, there exists a positive constant $C$
(independent of $l$) such that
$$
  \exp \left[ \int_{M \setminus \tilde{B}_i} \a G(x,y) h_l(y)
  dV_g(y) \right] \leq C, \qquad \hbox{ for any } x \in B_i.
$$
Then by \eqref{eq:expau} we have
\begin{equation}\label{eq:intAi}
    \int_{B_i} \exp \left[\a (v_l(x) - \ov{v}_l) \right] dV_g(x) \leq C \int_{B_i}
    \exp \left[ \int_{M} \a |G(x,y)| |h_l(y)| \chi_{\tilde{B}_i}
    dV_g(y) \right] dV_g(x).
\end{equation}
Now, as in \cite{bm}, we can use the Jensen's inequality to get
\begin{eqnarray*}
   \exp \left[ \int_M \a |G(x,y)| |h_l(y)| \chi_{\tilde{B}_i}
   dV_g(y) \right] \leq \int_M \exp \left[ \a \|h_l
   \chi_{\tilde{B}_i}\|_{L^1(M)} |G(x,y)| \right] \frac{|h_l
   \chi_{\tilde{B}_i}|(y)}{\|h_l \chi_{\tilde{B}_i}\|_{L^1(M)}}
  dV_g(y),
\end{eqnarray*}
and hence, by the Fubini Theorem and \eqref{eq:intAi}
$$
\int_{B_i} \exp \left[\a (v_l(x) - \ov{v}_l) \right] dV_g(x) \leq
C \sup_{y \in M} \int_M \exp \left[ \a \|h_l
\chi_{\tilde{B}_i}\|_{L^1(M)} |G(x,y)| \right] dV_g(x).
$$
By \eqref{eq:G}, there holds
$$
  \int_M \exp \left[ \a \|h_l \chi_{\tilde{B}_i}\|_{L^1(M)} |G(x,y)|
  \right] dV_g(x) \leq C \int_M \left( \frac{1}{|x-y|}
  \right)^{\frac{\a \|h_l \chi_{\tilde{B}_i}\|_{L^1(M)}}{8 \pi^2}} dV_g(x).
$$
The last integral is finite if
\begin{equation}\label{eq:intfin}
  \frac{\a \|h_l \chi_{\tilde{B}_i}\|_{L^1(M)}}{8 \pi^2} < 4 \qquad
  \Leftrightarrow \qquad \a \int_{\tilde{B}_i} |h_l| dV_g < 32 \pi^2.
\end{equation}
By \eqref{eq:concontr}, this is satisfied for some $\a > 4$
provided we take $l$ is sufficiently large. We have shown that
$\int_{B_i} e^{\a (v_l - \ov{v}_l)} dV_g < + \infty$ for every $i
= 1, \dots, L$. Since $M$ is covered finitely many $B_i$'s, the
conclusion follows.
\end{pf}

\begin{rem}\label{r:seqm}
Using the same proof, it is possible to extend Proposition
\ref{p:cc} to the case in which also the metric on $M$ depends on
$l$, and converges to some smooth $g$ in $C^m(M)$ for any integer
$m$. We have to use this variant in the next subsection.
\end{rem}

\subsection{Asymptotic profile}\label{ss:ap}

We consider now the alternative in Proposition \ref{p:cc} for
which compactness does not hold, applied to the case $h_l = 2 k_l
e^{4 u_l} - Q_l$. We assume that there exist $\rho \in \left( 0,
\frac{\pi^2}{k_0} \right)$, radii $(r_l)_l, (\hat{r}_l)_l$ and
points $(x_l)_l \subseteq M$ with the following properties
\begin{equation}\label{eq:trlxl}
    \hat{r}_l \to 0; \qquad \frac{r_l}{\hat{r}_l} \to 0;
    \qquad \int_{B_{r_l}(x_l)} e^{4 u_l} dV_g = \rho; \qquad
    \int_{B_{r_l}(y)} e^{4 u_l} dV_g < \frac{\pi^2}{k_0} \quad
    \hbox{ for every } y \in B_{\hat{r}_l}(x_l).
\end{equation}

\begin{rem}\label{r:ex}
If the second alternative in Proposition \ref{p:cc} holds, an
example of the situation described in \eqref{eq:trlxl} is the
following. Choose $r_l, x_l$ satisfying
\begin{equation}\label{eq:intxlrl}
    \int_{B_{r_l}(x_l)} e^{4 u_l} dV_g = \sup_{x \in M} \int_{B_{r_l}(x)}
    e^{4 u_l} dV_g = \rho.
\end{equation}
Then $r_l \to 0$ as $l \to + \infty$, and we can take $\hat{r}_l =
r_l^{\frac 12}$.
\end{rem}

\

\noindent Given a small $\d > 0$, we consider the exponential maps
$$
 \exp_l : B_\d^{\R^4} \to M; \qquad \qquad \exp_l(0) = x_l,
$$
where $B_\d^{\R^4} = \left\{ x \in \R^4 \; : \; |x| < \d
\right\}$. We also define the metric $\tilde{g}_l$ on
$B_\d^{\R^4}$ by $\tilde{g}_l := (\exp_l)^* g$, and the functions
$\tilde{u}_l : B_\d^{\R^4} \to \R$ by
$$
  \tilde{u}_l = u_l \circ \exp_l.
$$
Now in $\R^4$ we consider the dilation $T_l : x \mapsto r_l x$,
and we define another sequence
\begin{equation}\label{eq:hwl}
  \hat{u}_l(x) = \tilde{u}_l (T_l x) + \log r_l, \qquad \qquad x \in
  B_{\frac{\d}{r_l}}^{\R^4}.
\end{equation}
Using a change of variables, one easily verifies that the function
$\tilde{u}_l$ solves the equation
$$
P_{\tilde{g}_l} \tilde{u}_l(x) + 2 Q_{l} (x) = 2 k_l e^{4
\tilde{u}_l (x)}; \qquad \qquad x \in B_{\d}^{\R^4}.
$$
Hence, setting $\hat{g}_l = r_l^{-2} T_l^* \tilde{g}_l$ and using
the conformal properties of the Paneitz operator we obtain that
$\hat{u}_l$ satisfies
\begin{equation}\label{eq:eqhwl}
    P_{\hat{g}_l} \hat{u}_l(x) + 2 r_l^4 Q_l (T_l x) = 2
    k_l e^{4 \hat{u}_l(x)}; \qquad \qquad x \in
    B^{\R^4}_{\frac{\d}{r_l}}.
\end{equation}
Note that the metrics $\hat{g}_l$ converge in $C^m_{loc}(\R^4)$ to
the flat metric $(dx)^2$ for any integer $m$. Also, since
$(Q_l)_l$ are uniformly bounded functions on $M$, one also finds
$$
  r_l^4 Q_{\tilde{g}_l} (T_l \cdot) \to 0 \qquad \qquad \hbox{ in }
  C^0_{loc}(\R^4).
$$
By \eqref{eq:trlxl}, using a change of variables we obtain
\begin{equation}\label{eq:inthwl}
    \rho = \int_{B_{r_l}(x_l)} e^{4 u_l} dV_g =
    \int_{\frac{1}{r_l} (\exp_l)^{-1} B_{r_l}(x_l)} e^{4
    \hat{u}_l} d V_{\hat{g}_l},
\end{equation}
where $o_l(1) \to 0$ as $l \to + \infty$. Note also that the sets
$\frac{1}{r_l} (\exp_l)^{-1} B_{r_l}(x_l) \subseteq \R^4$ approach
the unit ball $B_1^{\R^4}$ as $l \to + \infty$. Moreover, by the
last inequality in \eqref{eq:trlxl} and by our choice of $\rho$,
it is easy to derive that
\begin{equation}\label{eq:nonc3}
  \int_{B_{\frac 12}^{\R^4}(y)} e^{4 \hat{u}_l} d V_{\hat{g}_l} <
  \frac{\pi^2}{k_0} \qquad \qquad \hbox{ for every } y \in
  B^{\R^4}_{\frac{\hat{r}_l}{2 r_l}}.
\end{equation}

\

\noindent Regarding the functions $\hat{u}_l$, we have the
following convergence result.

\begin{pro}\label{p:hwlb}
Suppose $\rho \in \left( 0, \frac{\pi^2}{k_0} \right)$, $(r_l)_l,
(\tilde{r}_l)_l$, $(x_l)_l$ and $(u_l)_l$ satisfy
\eqref{eq:trlxl}, and let $(\hat{u}_l)_l$ be defined by
\eqref{eq:hwl}. Then there exists $\l > 0$, $x_0 \in \R^4$ and $\a
\in (0,1)$ such that
$$
  \hat{u}_l \to \hat{u}_\infty \qquad \qquad \hbox{ in }
  C^\a_{loc}(\R^4) \hbox{ and in } H^2_{loc}(\R^4)
$$
for some $\a \in (0,1)$, where the function $\hat{u}_\infty$ is
given by
\begin{equation}\label{eq:hwinf}
    \hat{u}_\infty(x) = \log \frac{2 \l}{1 + \l^2|x - x_0|^2} -
  \frac 14 \log \left( \frac 13 k_0 \right); \qquad x \in \R^4.
\end{equation}
Moreover, if $b_l \to + \infty$ sufficiently slowly, one has
\begin{equation}\label{eq:intalrlu}
    \int_{B_{b_l r_l}(x_l)} e^{4 u_l} d V_g \to \frac{8 \pi^2}{k_0}
  \qquad \hbox{ as } l \to + \infty.
\end{equation}
\end{pro}

\begin{pf}
Given $R > 0$, we define a smooth cut-off function $\Psi_R$
satisfying
$$
  \left\{%
\begin{array}{ll}
    \Psi_R(x) = 1, & \hbox{ for } |x| \leq \frac R2;  \\
    \Psi_R(x) = 0, & \hbox{ for } |x| \geq R.\\
\end{array}%
\right.
$$
We also set
$$
 a_l = \frac{1}{|B_R^{\R^4}|} \int_{B_R^{\R^4}} \hat{u}_l d V_{\hat{g}_l}; \qquad \qquad v_l
 = \Psi_R \hat{u}_l + (1 - \Psi_R) a_l = a_l + \Psi_R (\hat{u}_l - a_l);
$$
$$
  \hat{v}_l = v_l - a_l.
$$
We notice that the functions $v_l$ coincide with $a_l$ outside
$B_R^{\R^4}$ and that $\hat{v}_l$ is identically zero outside
$B_R^{\R^4}$. By Lemma \ref{l:lodlp} and some scaling argument one
finds
\begin{equation}\label{eq:intp1}
  \int_{B^{\R^4}_{2R}} \left( |\n^3 \hat{u}_l|^p + |\n^2 \hat{u}_l|^p
  + |\n \hat{u}_l|^p  \right) dV_{\hat{g}_l} \leq C_R; \qquad l
   \in \N, p \in \left( 1, \frac 43 \right),
\end{equation}
and hence by the Poincar\'e inequality (recall that the
$\hat{v}_l$'s have a uniform compact support) it follows that
\begin{equation}\label{eq:intp2}
  \int_{B^{\R^4}_{R}} |\hat{v}_l|^p dV_{\hat{g}_l} \leq C_R; \qquad
  \qquad l \in \N, p \in \left( 1, \frac 43 \right).
\end{equation}
By \eqref{eq:eqhwl} there holds
\begin{eqnarray}\label{eq:intp22}
    P_{\hat{g}_l} \hat{v}_l & = & (\D_{\hat{g}_l})^2 [\Psi_R (\hat{u}_l -
    a_l)] + L_l [\Psi_R (\hat{u}_l - a_l)] = \Psi_R P_{\hat{g}_l}
    \hat{u}_l + \tilde{L}_l (\hat{u}_l - a_l) \nonumber \\ & = & 2 k_l
    \Psi_R e^{4 \hat{u}_l} + \hat{f}_l,
\end{eqnarray}
where
$$
  \hat{f}_l = \tilde{L}_l (\hat{u}_l - a_l) - 2 r_l^4
  Q_l (T_l \cdot).
$$
Here $(L_l)_l$ are linear operators which contain derivatives of
order $1$ and $2$ with uniformly bounded and smooth coefficients.
Also, $(\tilde{L}_l)_l$ are linear operators which contain
derivatives of order $0, 1, 2$ and $3$ with uniformly bounded and
smooth coefficients. As a consequence, by \eqref{eq:intp1} and
\eqref{eq:intp2} one has
\begin{equation}\label{eq:intp3}
  \int_{B_{2R}^{\R^4}} |\hat{f}_l|^p dV_{\hat{g}_l} \leq C_R;
  \qquad l \in \N, p \in \left( 1, \frac 43 \right).
\end{equation}
Hence using \eqref{eq:nonc3} and Remark \ref{r:seqm} one finds
\begin{equation}\label{eq:intqhvl}
    \int_{B_R^{\R^4}} e^{4 q \hat{v}_l} dV_{\hat{g}_l} \leq C
    \qquad \qquad \hbox{ for some } q > 1
\end{equation}
and for some fixed constant $C$. Remark \ref{r:seqm} applies
indeed to the case of a compact manifold while in the present
situation we are working in $\R^4$ (endowed with the metric
$\hat{g}_l$). But since all the functions $\hat{v}_{l}$ vanish
identically outside $B_R^{\R^4}$, we can embed a fixed
neighborhood of $(B_{2R}^{\R^4}, \hat{g}_l)$ into a compact
manifold, a torus for example, such that its metric (coinciding
with $\hat{g}_l$ on $B_{2R}^{\R^4}$) converges to the flat one.

On the other hand, from \eqref{eq:nonc3} we deduce
$$
   a_l = \frac{1}{|B_R^{\R^4}|} \int_{B_R^{\R^4}} \hat{u}_l dV_{\hat{g}_l} \leq
   \frac{1}{4 |B_R^{\R^4}|} \int_{B_R^{\R^4}} e^{4 \hat{u}_l} dV_{\hat{g}_l} \leq C,
$$
and from \eqref{eq:inthwl}, since $v_l = \hat{u}_l$ in
$B_R^{\R^4}$
$$
C^{-1} \leq \int_{B_R^{\R^4}} e^{4 v_l} dV_{\hat{g}_l} \leq e^{4
a_l} \int_{B_R^{\R^4}} e^{4 \hat{v}_l} dV_{\hat{g}_l} \leq C e^{4
a_l}.
$$
This implies $a_l \geq - C$, and hence we find
$$
|a_l| \leq C.
$$
As a consequence of this estimate and \eqref{eq:intqhvl} we get
the following uniform improved integrability for $\hat{u}_l$
(recall the definition of $v_l$ and $\hat{v}_l$)
$$
  \int_{B_R^{\R^4}} e^{4 q \hat{u}_l} dV_{\hat{g}_l} \leq C, \qquad \hbox{
  for some } q > 1.
$$
This estimate, joint with \eqref{eq:intp1}, \eqref{eq:intp22},
\eqref{eq:intp3} and standard elliptic regularity results, yields
that $\hat{u}_l$ is bounded in $W^{4,q}(B_{\frac R2}^{\R^4})$.
Hence, by the arbitrarity of $R$, $(\hat{u}_l)_l$ converge
strongly in $C^\a_{loc}(\R^4)$ for some $\a \in (0,1)$ and
strongly in $H^2_{loc}(\R^4)$ to a function $\hat{u}_\infty \in
C^\a_{loc}(\R^4) \cup H^2_{loc}(\R^4)$.

Now we prove that $\hat{u}_\infty$ has the form in
\eqref{eq:hwinf}. First of all, we test equation \eqref{eq:eqhwl}
on a smooth function $\var$ with compact support. Integrating by
parts we obtain
$$
  \langle P_{\hat{g}_l} \hat{u}_l, \var \rangle + 2 r_l^4 \int_{\R^4}
  Q_{l} (T_l \cdot) \var d V_{\hat{g}_l} = 2 k_l \int_{\R^4}
  e^{4 \hat{u}_l} \var d V_{\hat{g}_l}.
$$
As $l$ tends to infinity we get
$$
  \langle P_{\R^4} \hat{u}_\infty, \var \rangle = 2 k_0
  \int_{\R^4} e^{4 \hat{u}_\infty} \var d V_{\R^4} + o_l(1).
$$
Hence the limit function $\hat{u}_\infty$ satisfies
\begin{equation}\label{eq:r4}
    \D^2_{\R^4} \hat{u}_\infty = 2 k_0 e^{4 \hat{u}_\infty}
    \qquad \qquad \hbox{ in } \R^4,
\end{equation}
and, by semicontinuity
\begin{equation}\label{eq:cond}
    \int_{\R^4} e^{4 \hat{u}_\infty} dV_{\R^4} \leq 1,
\end{equation}
since by \eqref{eq:noul} and some scaling there holds
$\int_{B^{\R^4}_{\frac{\d}{r_l}}} e^{4 \hat{u}_l} dV_{\hat{g}_l}
\leq 1$.

The solutions of \eqref{eq:r4}-\eqref{eq:cond}, with $k_0 > 0$,
have been classified in \cite{lin}, and one of the following two
possibilities occur

\

{\bf (a)} $\hat{u}_\infty$ is of the form \eqref{eq:hwinf}, or

\

{\bf (b)} $\D_{\R^4} \hat{u}_\infty$ has the following asymptotic
behavior
\begin{equation}\label{eq:hunbub}
    - \D_{\R^4} \hat{u}_\infty(x) \to a > 0 \qquad \hbox{ for }
    |x| \to + \infty.
\end{equation}
Following \cite{rs}, we show that the second alternative does not
happen. In fact, assuming {\bf (b)}, for $R$ large we have
\begin{equation}\label{eq:BL}
    \lim_{l \to + \infty} \int_{B_R^{\R^4}} (- \D_{\hat{g}_l}
    \hat{u}_l) dV_{\hat{g}_l} = \int_{B_R^{\R^4}} (- \D_{\R^4}
    \hat{u}_\infty) dV_{\R^4} \sim \frac{\o_3}{4} a R^4,
\end{equation}
where $\o_3 = |S^3| = 2 \pi^2$. Scaling back to $M$ (recall that
the dilation factor is $r_l$), we obtain
\begin{equation}\label{eq:BL2}
    \lim_{l \to + \infty} \int_{B_{R r_l}(x_l)} (- \D u_l) dV_g
    \sim \ov{C} a R^4 r_l^2,
\end{equation}
for some $\ov{C} > 0$. On the other hand, by Lemma \ref{l:lodlp}
we get
\begin{equation}\label{eq:c2}
  \int_{B_{Rr_l}(x_l)} (- \D u_l) dV_g \leq \hat{C}_0 r_l^2
  R^2.
\end{equation}
Taking $R$ sufficiently large, from \eqref{eq:BL2} and
\eqref{eq:c2} we reach a contradiction.

Hence the alternative {\bf (a)} holds and $\hat{u}_\infty$ arises
as a conformal factor of a stereographic projection of $S^4$ onto
$\R^4$, which must satisfy
\begin{equation}\label{eq:intehwi}
    \int_{\R^4} e^{4 \hat{u}_\infty} dV_{\R^4} = \frac{8 \pi^2}{k_0}.
\end{equation}
This concludes the proof.
\end{pf}

\section{Simple blow-ups}\label{s:sbu}

In this section we consider an unbounded sequence of solutions
$(u_l)_l$ and we examine a particular class of blow-up points,
which we call {\em simple}, in analogy with a definition
introduced by R.Schoen. In Proposition \ref{p:sbu} below we give
some quantitative estimate on the concentration at simple blow-up
points. Then in the next section we show that every general
blow-up phenomenon can be essentially reduced to the study of
finitely many simple blow-ups. In the following $i(M)$ denotes the
injectivity radius of $M$.

\begin{df}\label{d:sbu}
If $(u_l)_l$ satisfies \eqref{eq:pl} and \eqref{eq:noul}, we say
that the three sequences $(x_l)_l \subseteq M$, $r_l \to 0$,
$(s_l)_l \subseteq \R_+$, $|s_l| \leq i(M)$ are a {\em simple
blow-up} for $(u_l)_l$ if the following properties hold
\begin{equation}\label{eq:rlslxl}
    \frac{s_l}{r_l} \to + \infty; \qquad \quad \exists R_l \to + \infty
    \hbox{ s.t. } \left\| \hat{u}_l
    \to \log \frac{2}{1 + |\cdot|^2} - \frac{1}{4} \log \left( \frac 13 k_0
    \right) \right\|_{H^4(B_{R_l}^{\R^4}) \cap
    C^\a(B_{R_l}^{\R^4})} \to 0;
\end{equation}
\begin{equation}\label{eq:al}
    \forall \rho > 0 \exists C_\rho > 0 \hbox{ s.t. if } \int_{B_s(y)} e^{4 u_l}
    d V_g \geq \rho \hbox{ with } B_s(y) \subseteq
    B_{s_l}(x_l) \setminus B_{R_l r_l}(x), \hbox{ then } s \geq
    C_\rho^{-1} |y - x_l|,
\end{equation}
where $\hat{u}_l$ is defined in \eqref{eq:hwl}.
\end{df}

\

\noindent The main result of this section is the following
proposition.

\begin{pro}\label{p:sbu} Suppose $(x_l)_l$, $(r_l)_l$, $(s_l)_l$
are a simple blow-up for $(u_l)_l$. Then there exists a fixed $C >
0$ such that
\begin{equation}\label{eq:estsbu}
    \int_{B_{C^{-1} s_l}(x_l)} e^{4
    u_l} dV_g = \frac{8 \pi^2}{k_0} + o_l(1),
\end{equation}
where $o_l(1) \to 0$ as $l \to + \infty$.
\end{pro}

\

\begin{rem}\label{r:sbu} (a) We notice that, if $\hat{u}_l$
satisfies the assertion in Proposition \ref{p:hwlb}, it is always
possible to modify $(x_l)_l$ and $(r_l)_l$ in order to get $x_0 =
0$ and $\l = 1$.

(b) Proposition \ref{p:sbu} is basically an improvement of formula
\eqref{eq:al} to a sequence of sets with larger size.
\end{rem}

\

\noindent The proof of Proposition \ref{p:sbu} is based on the
analysis of the next two subsections. In the first we prove some
Harnack inequality in integral form while in the second, defining
\begin{equation}\label{eq:}
    A_{r,l} = \left\{ x \in M \; : \; r < |x - x_l| < 2 r
    \right\},
\end{equation}
we study the average of $u_l$ on $A_{r,l}$ as a function of $r$.

\subsection{Integral Harnack-type inequalities}\label{ss:ih}

In this subsection we prove some integral Harnack-type
inequalities for the functions $(u_l)_l$ near simple blow-ups.
Although it is maybe possible to get pointwise estimates on the
solutions, for our purposes it is sufficient to get integral
volume estimates. We need first a preliminary result involving the
average of the Green's function $G$ on annuli. Given $\rho \in
\left( 0, \frac{\pi^2}{k_0} \right)$, let $C_\rho$ be the
corresponding constant in \eqref{eq:al} (which we can suppose
bigger than $1$), and we define the following sets
\begin{equation}\label{eq:2}
    A'_{r,l} = \left\{ x \in M \; : \; \frac{5}{4} r < |x - x_l| <
    \frac{7}{4} r \right\} \subseteq A_{r,l}; \qquad r \in
    \left( a_l r_l, s_l \right);
\end{equation}
\begin{equation}\label{eq:btb}
    \mathcal{B}_{r}(x) = B_{\frac{r}{16 C_\rho}}(x) \subseteq A'_{r,l};
    \qquad \tilde{\mathcal{B}}_{r}(x) = B_{\frac{r}{8 C_\rho}}(x)
    \subseteq A'_{r,l}, \qquad \qquad x \in A'_{r,l}.
\end{equation}

\begin{lem}\label{l:Green}
Suppose $(x_l)_l \subseteq M$, $(s_l)_l \subseteq \R_+$, $|s_l|
\leq i(M)$, and let $A_{r,l}, A'_{r,l}$,
$\tilde{\mathcal{B}}_{r}(x)$ be defined respectively in
\eqref{eq:}, \eqref{eq:2} and \eqref{eq:btb}. Then there exists a
positive constant $C = C(C_\rho)$, independent of $r$ and $l$ such
that, setting
$$
  f_{r,l}(y) = \frac{1}{|A_{r,l}|} \int_{A_{r,l}} G(z,y) dV_g(z),
$$
there holds
\begin{equation}\label{eq:estfry}
    \left\{%
\begin{array}{ll}
    \left| f_{r,l}(y) - \frac{1}{8 \pi^2} \log \frac 1r \right| \leq C
    & \hbox{ for every } x \in A'_{r,l}, y \in
    \tilde{\mathcal{B}}_{r}(x),
    \\ \left| f_{r,l}(y) - G(x,y) \right| \leq C & \hbox{ for every }
    x \in A'_{r,l}, y \in M \setminus \tilde{\mathcal{B}}_{r}(x); \\
\end{array}%
\right. \qquad r \leq i(M).
\end{equation}
\end{lem}

\begin{pf}
We first notice that the following inequality holds
\begin{equation}\label{eq:logeucl}
    \left| \ov{f}_r(y) - \log \frac 1r \right| \leq \ov{C}; \qquad
    \qquad |y| \leq 4 r,
\end{equation}
where
$$
  A_r = \{ x \in \R^4 \; : r < |x| < 2r \}; \qquad \qquad \ov{f}_r(y) =
  \frac{1}{|A_r|_{\R^4}} \int_{A_r} \log \frac{1}{|z - y|_{\R^4}}
  dV_{\R^4}.
$$
Here $|A_r|_{\R^4}$ stands for the Lebesgue measure of $A_r$ and
$|z - y|_{\R^4}$ denotes the Euclidean distance.

The inequality is indeed trivial for $r = 1$ since $\ov{f}_1(y)$
is bounded on $B_4^{\R^4}$, while for a general $r$ it is
sufficient to use a scaling argument. We use \eqref{eq:G}, the
exponential map and standard geometric estimates on $M$ (see
\eqref{eq:dVg} below for the volume element) to write
\begin{eqnarray*}
    8 \pi^2 f_{r,l}(y) & = & \frac{1}{|A_{r,l}|} \int_{A_{r,l}}
    \log \frac{1}{|y - z|} dV_g(z) + O(1) \\ & = & (1 + O(r^2))
    \frac{1}{|A_r|_{\R^4}}
    \int_{A_r} \log \frac{1}{|y - z|_{\R^4}} (1 + O(r^2)) dV_{\R^4} +
    O(1) \\ & = & (1 + O(r^2)) \ov{f}_r(y) + O(1); \qquad \qquad
    \qquad y \in B_{4r}(x_l).
\end{eqnarray*}
Jointly with \eqref{eq:logeucl}, this proves the first estimate in
\eqref{eq:estfry}.

The second one is trivial for $y \in B_{4r}(x_l) \setminus
\tilde{\mathcal{B}}_{r}(x)$, by the preceding argument. For $y \in
M \setminus B_{4r}(x_l)$, we notice that
$$
  C^{-1} \leq \frac{|z - y|}{|x - y|} \leq C \qquad \qquad \hbox{
  for } z \in A_{r,l}, x \in A'_{r,l},
$$
and we use again \eqref{eq:G}. This concludes the proof.
\end{pf}

\

\noindent Next, we prove some inequality involving the integral of
the function $e^{4 u_l}$ and the average of $u_l$ on small annuli.
Recall the definitions of $A_{r,l}$ and $A'_{r,l}$ in \eqref{eq:}
and \eqref{eq:2}, and those of $\mathcal{B}_r(x),
\tilde{\mathcal{B}}_r(x)$ in \eqref{eq:btb}.

\begin{lem}\label{l:trueintw}
Suppose that $(x_l)_l \subseteq M$, $r_l \to 0$, $(s_l)_l
\subseteq \R_+$, $|s_l| \leq i(M)$ are a simple blow-up for
$(u_l)_l$. Suppose $a_l \to + \infty$, and define
$$
   \ov{u}_{r,l} = \frac{1}{|A_{r,l}|} \int_{A_{r,l}} u_l
   dV_g; \qquad \qquad a_l r_l < r < s_l.
$$
Then, if $l$ is sufficiently large, there exists a positive
constant $C$ (independent of $l$ and $r$) such that
$$
  \int_{A'_{r,l}} e^{4 u_l} dV_g \leq C |A_{r,l}| e^{4
  \ov{u}_{l,r}}; \qquad \qquad a_l r_l < r < s_l.
$$
\end{lem}

\begin{pf} Using \eqref{eq:Grepwl} and recalling the definition of
$f_l$ (see \eqref{eq:fl}) and that of $f_{r,l}$ (see Lemma
\ref{l:Green}), we have
$$
  \ov{u}_{r,l} = \ov{u}_l + \int_M f_{r,l}(y) f_l(y) dV_g(y).
$$
For $x \in A'_{r,l}$, we divide the last integral into
$\tilde{\mathcal{B}}_{r}(x)$ and its complement, to obtain
$$
  \exp \left( 4 (\ov{u}_{r,l} - \ov{u}_l) \right) = \exp \left( 4
  \int_{\tilde{\mathcal{B}}_{r}(x)} f_{r,l}(y) f_l(y) dV_g(y)
  \right) \exp \left( 4 \int_{M \setminus \tilde{\mathcal{B}}_{r}(x)}
  f_{r,l}(y) f_l(y) dV_g(y) \right).
$$
Using Lemma \ref{l:Green} and the fact that $(f_l)_l$ is bounded
in $L^1(M)$, we then find
$$
  \exp \left( 4 (\ov{u}_{r,l} - \ov{u}_l) \right) \geq C^{-1} \exp
  \left( \frac{1}{2 \pi^2} \log \frac 1r \int_{\tilde{\mathcal{B}}_{r}(x)}
  f_l(y) dV_g(y) \right)  \exp \left( 4
  \int_{M \setminus \tilde{\mathcal{B}}_{r}(x)} G(x,y) f_l(y) dV_g(y) \right).
$$
Hence, integrating on $A_{r,l}$ we obtain
\begin{eqnarray}\label{eq:intavwl}
    \int_{A_{r,l}} e^{4 (\ov{u}_{r,l} - \ov{u}_l)} dV_g \geq
    C^{-1} |A_{r,l}| \left( \frac 1r \right)^{\frac{\int_{\tilde{\mathcal{B}}_{r}(x)}
  f_l dV_g}{2 \pi^2}} \exp \left( 4
  \int_{M \setminus \tilde{\mathcal{B}}_{r}(x)} G(x,y) f_l(y) dV_g(y)
  \right).
\end{eqnarray}
On the other hand, again by \eqref{eq:Grepwl}, for $x \in
A'_{r,l}$ and a.e. $z \in \mathcal{B}_r(x)$ we have also
$$
  u_l(z) - \ov{u}_l = \int_{M \setminus \tilde{\mathcal{B}}_{r}(x)}
  G(z,y) f_l(y) dV_g(y) + \int_{\tilde{\mathcal{B}}_{r}(x)} G(z,y) f_l(y) dV_g(y).
$$
Then, exponentiating and integrating on $\mathcal{B}_r(x)$ we get
\begin{eqnarray}\label{eq:1+2}
    & & \int_{\mathcal{B}_r(x)} e^{4 (u_l(z) - \ov{u}_l)} dV_g(z) \\ & = &
    \int_{\mathcal{B}_r(x)} \exp \left( 4
    \int_{M \setminus \tilde{\mathcal{B}}_{r}(x)}
  G(z,y) f_l(y) dV_g(y) \right) \exp \left( 4 \int_{\tilde{\mathcal{B}}_{r}(x)}
  G(z,y) f_l(y) dV_g(y) \right) dV_g(z) \nonumber \\ & \leq &
  \underbrace{\sup_{z \in \mathcal{B}_r(x)} \exp \left( 4
    \int_{M \setminus \tilde{\mathcal{B}}_{r}(x)} G(z,y) f_l(y) dV_g(y)
    \right)}_{J}
    \underbrace{ \int_{\mathcal{B}_r(x)} \exp
    \left( 4 \int_{\tilde{\mathcal{B}}_{r}(x)} G(z,y) f_l(y) dV_g(y) \right)
    dV_g(z)}_{JJ}. \nonumber
  \nonumber
\end{eqnarray}
Now we write
$$
  \int_{M \setminus \tilde{\mathcal{B}}_{r}(x)}
  G(z,y) f_l(y) dV_g(y) = \int_{M \setminus \tilde{\mathcal{B}}_{r}(x)}
  G(x,y) f_l(y) dV_g(y) + \int_{M \setminus \tilde{\mathcal{B}}_{r}(x)}
  \left( G(z,y) - G(x,y) \right) f_l(y) dV_g(y).
$$
Using \eqref{eq:G}, for $z \in \mathcal{B}_r(x)$ and $y \in M
\setminus \tilde{\mathcal{B}}_{r}(x)$, we have
$$
  G(z,y) - G(x,y) = O(1) + \frac{1}{8 \pi^2}
  \log \frac{|z - y|}{|x - y|} = O(1).
$$
As a consequence we deduce
\begin{equation}\label{eq:1+1}
  J \leq C \exp \left( 4 \int_{M \setminus
  \tilde{\mathcal{B}}_{r}(x)} G(x,y) f_l(y) dV_g(y) \right).
\end{equation}
We now turn to $JJ$. Since $z \in \mathcal{B}_r(x)$ and $y \in
\tilde{\mathcal{B}}_{r}(x)$, $G(z,y)$ is positive (for $r$
sufficiently small), and hence
$$
  \int_{\tilde{\mathcal{B}}_{r}(x)} G(z,y) f_l(y) dV_g(y) \leq
  \int_{\tilde{\mathcal{B}}_{r}(x)} G(z,y) |f_l|(y) dV_g(y).
$$
Using the Jensen inequality, as in the proof of Proposition
\ref{p:cc}, we obtain
\begin{eqnarray*}
    \exp \left( 4 \int_{\tilde{\mathcal{B}}_{r}(x)}
  G(z,y) f_l(y) dV_g(y) \right)
    \leq \int_{\tilde{\mathcal{B}}_{r}(x)} \exp \left(
    4 G(z,y) \|f_l\|_{L^1(\tilde{\mathcal{B}}_{r}(x))}
    \right)
    \frac{|f_l(y)|}{\|f_l\|_{L^1(\tilde{\mathcal{B}}_{r}(x))}} dV_g(y).
\end{eqnarray*}
Again \eqref{eq:G} implies
\begin{eqnarray*}
    JJ
    & \leq & \int_{\mathcal{B}_r(x)} dV_g(z) \int_{\tilde{\mathcal{B}}_{r}(x)}
    \exp \left( 4 G(z,y) \|f_l\|_{L^1(\tilde{\mathcal{B}}_{r}(x))}
    \right) \frac{|f_l(y)|}{\|f_l\|_{L^1(\tilde{\mathcal{B}}_{r}(x))}}
    dV_g(y) \\ & \leq & C \int_{\mathcal{B}_r(x)} dV_g(z)
    \int_{\tilde{\mathcal{B}}_{r}(x)} \left( \frac{1}{|z - y|}
    \right)^{\frac{\|f_l\|_{L^1(\tilde{\mathcal{B}}_{r}(x))}}{2 \pi^2}}
    \frac{|f_l(y)|}{\|f_l\|_{L^1(\tilde{\mathcal{B}}_{r}(x))}}
    dV_g(y).
\end{eqnarray*}
Now, the Fubini theorem and some elementary computations yield
\begin{equation}\label{eq:2+2}
    JJ \leq C \sup_{y \in M} \int_{\mathcal{B}_r(x)} dV_g(z) \left( \frac{1}{|z - y|}
    \right)^{\frac{\|f_l\|_{L^1(\tilde{\mathcal{B}}_{r}(x))}}{2
    \pi^2}} \leq C r^{4 - \frac{\|f_l\|_{L^1(\tilde{\mathcal{B}}_{r}(x))}}{2
    \pi^2}}.
\end{equation}
In the last inequality we have used the fact that
$\|f_l\|_{L^1(\tilde{\mathcal{B}}_{r}(x))}$ is uniformly small
since we are dealing with a simple blow-up, see \eqref{eq:al}, and
since we have chosen $\tilde{\mathcal{B}}_{r}(x)$ suitably. This
implies that the last constant $C$ is independent of $r$ and $l$.
From \eqref{eq:1+2}, \eqref{eq:1+1} and \eqref{eq:2+2} it follows
that
$$
\int_{\mathcal{B}_r(x)} e^{4 (u_l(z) - \ov{u}_l)} dV_g(z) \leq
 C r^{4 - \frac{\|f_l\|_{L^1(\tilde{\mathcal{B}}_{r}(x))}}{2
    \pi^2}} \exp \left( 4 \int_{M \setminus \tilde{\mathcal{B}}_{r}(x)}
  G(x,y) f_l(y) dV_g(y) \right).
$$
Now the assertion of the Lemma follows from the last formula,
\eqref{eq:intavwl} and the observation that, since $f_l = 2k_l
e^{4u_l} - 2 Q_l$, it is
$\|f_l\|_{L^1(\tilde{\mathcal{B}}_{r}(x))} =
\int_{\tilde{\mathcal{B}}_{r}(x)} f_l dV_g + O(r^4)$, and hence
$$
  r^{4 - \frac{\|f_l\|_{L^1(\tilde{\mathcal{B}}_{r}(x))}}{2 \pi^2}} \leq
  C |A_{r,l}| \left( \frac 1r \right)^{\frac{\int_{\tilde{\mathcal{B}}_{r}(x)}
  f_l dV_g}{2 \pi^2}} \qquad \hbox{ independently of } r \hbox{ and } l.
$$
This concludes the proof.
\end{pf}

\

\noindent Next we show some further estimates involving the
Laplacian of $u_l$. Recall that we have set $f_l = 2 k_l e^{4 u_l}
- Q_l$, see \eqref{eq:fl}.

\begin{lem}\label{l:harlap}
Suppose that $(x_l)_l \subseteq M$, $(\Sig_l)_l, (S_l)_l \subseteq
\R_+$, $i(M) \geq S_l > \Sig_l > 0$, and that $(u_l)_l$ satisfies
\eqref{eq:pl} and \eqref{eq:noul}. Suppose also that
$$
  \int_{B_{S_l}(x_l) \setminus B_{\Sig_l}(x_l)} e^{4 u_l} d V_g \leq
  \e.
$$
Then, for any $R > 0$ sufficiently large and any $r \in [\Sig_l +
R, S_l - R]$, one has
$$
  \int_{A_{r,l}} |x-x_l|^2 (- \D u_l(x)) dV_g(x) = \left( \frac{15}{8}
  \int_{B_{\frac rR}(x_l)} f_l dV_g + o_R(1) + O(\e R^2) + o_r(1) \right)
  r^4.
$$
where $o_R(1) \to 0$ as $R \to + \infty$ and $o_r(1) \to 0$ as $r
\to 0$.
\end{lem}

\begin{pf}
We can write \eqref{eq:pl} in the following form
$$
  -\D (- \D u_l) = f_l + F_l(u_l),
$$
where $F_l$ is a linear expression in $\n u_l$ and $\n^2 u_l$ with
uniformly bounded coefficients. If $\hat{G}$ is the Green's
function for the (negative) laplacian on $M$, then it is a
standard fact that
\begin{equation}\label{eq:hatG}
    \hat{G}(x,y) = (1 + o(1))\frac{1}{4 \pi^2 |x-y|^2};
    \qquad (x, y) \in M \times M \setminus \hbox{diag},
\end{equation}
where $o(1) \to 0$ as $|x - y| \to 0$, see for example \cite{aul}.
Hence, using the representation formula, for a.e. $x \in A_{r,l}$
we obtain
\begin{equation}\label{eq:dwlv1v2}
  - \D u_l(x) = \int_M \hat{G}(x,y) f_l(y) dV_g(y) +
  \int_M \hat{G}(x,y) F_l(u_l)(y) dV_g(y) := v_{1,l}(x) + v_{2,l}(x).
\end{equation}
Given $R > 0$ large but fixed and for $|x - x_l| = r \in [\Sig_l +
R, S_l - R]$, we write
\begin{eqnarray*}
    v_{1,l}(x) & = & \int_{B_{\frac rR}(x_l)} \hat{G}(x,y) f_l(y) dV_g(y)
  + \int_{B_{Rr}(x_l) \setminus B_{\frac rR}(x_l)} \hat{G}(x,y) f_l(y)
  dV_g(y) \\ & + & \int_{M \setminus B_{Rr}(x_l)} \hat{G}(x,y) f_l(y)
  dV_g(y).
\end{eqnarray*}
From the asymptotics in \eqref{eq:hatG} and some scaling argument
we obtain (for $x \in A_{r,l}$)
$$
  \int_{B_{\frac rR}(x_l)} \hat{G}(x,y) f_l(y) dV_g(y) = (1 +
  o_r(1) + o_R(1))\frac{1}{4 \pi^2 r^2}
  \int_{B_{\frac rR}} f_l dV_g;
$$
$$
   \left| \int_{M \setminus B_{Rr}(x_l)}
  \hat{G}(x,y) f_l(y) dV_g(y) \right| \leq \frac{C}{(R r)^2},
$$
where $o_r(1) \to 0$ as $r \to 0$ and $o_R(1) \to 0$ as $R \to +
\infty$. Moreover, by our assumptions and \eqref{eq:fl}, we have
$$
  \int_{B_{Rr}(x_l) \setminus B_{\frac rR}(x_l)}
  f_l(y) dV_g(y) \leq C \e; \qquad \qquad f_l(x) \geq - C,
$$
where $C$ is independent of $r$, and $l$. Using the Fubini theorem
and reasoning as in the proof of Lemma \ref{l:lodlp} it follows
that
$$
  \left| \int_{A_{r,l}} dV_g(x) \int_{B_{Rr}(x_l) \setminus
  B_{\frac rR}(x_l)} \hat{G}(x,y) f_l(y) dV_g(y) \right| \leq C
  \e R^2 r^2.
$$
The last formulas imply
\begin{eqnarray}\label{eq:estv1l} \nonumber
  \int_{A_{r,l}} |x - x_l|^2 v_{1,l}(x) dV_g(x) & = & \left(
  \frac{1 + o_r(1) + o_R(1)}{4 \pi^2} \int_{B_{\frac rR}} f_l dV_g
  + O(\e R^2) + O \left( \frac{1}{R^2} \right) \right) |A_{r,l}|
  \\ & = & \left( \frac{15}{8} \int_{B_{\frac rR}(x_l)} f_l dV_g
  + o_R(1) + O(\e R^2) + o_r(1) \right) r^4.
\end{eqnarray}
To study the integral of $v_{2,l}$, we use again the
representation formula and we write
\begin{eqnarray*}
    |v_{2,l}(x)| & \leq & C \int_{B_{r^2}(x_l)} \frac{1}{|x-y|^2}
  \left( |\n^2 u_l|(y) + |\n u_l|(y) \right) dV_g(y) \\ & + &
  C \int_{M \setminus B_{\sqrt{r}}(x_l)} \frac{1}{|x-y|^2}
  \left( |\n^2 u_l|(y) + |\n u_l|(y) \right) dV_g(y) \\ & + &
  C \underbrace{\int_{B_{\sqrt{r}}(x_l) \setminus B_{r^2}(x_l)}
  \frac{1}{|x-y|^2} \left( |\n^2 u_l|(y) + |\n u_l|(y) \right)
  dV_g(y)}_{JJJ}.
\end{eqnarray*}
To estimate the first and the second integral, we notice that $|x
- y| \geq C^{-1} r$ and $|x - y| \geq C^{-1} \sqrt{r}$ for
respectively $y \in B_{r^2}(x_l)$ and $y \in B_{\sqrt{r}}(x_l)$
(recall that $x \in A_{r,l}$). Hence using Lemma \ref{l:lodlp} it
follows that
$$
  \int_{B_{r^2}(x_l)} \frac{1}{|x-y|^2}
  \left( |\n^2 u_l|(y) + |\n u_l|(y) \right) dV_g(y) \leq C r^2;
$$
$$
  \int_{M \setminus B_{\sqrt{r}}(x_l)} \frac{1}{|x-y|^2}
  \left( |\n^2 u_l|(y) + |\n u_l|(y) \right) dV_g(y) \leq \frac
  Cr.
$$
To estimate the third integral we use the H\"older's inequality to
find, for $\frac 1p + \frac{1}{p'} = 1$
$$
  JJJ \leq C \left( \int_{B_{\sqrt{r}}(x_l) \setminus B_{r^2}(x_l)}
  \frac{1}{|x-y|^{2p}} dV_g(y) \right)^{\frac 1p} \left(
  \int_{B_{\sqrt{r}}(x_l) \setminus B_{r^2}(x_l)} \left(
  |\n^2 u_l|(y) + |\n u_l|(y) \right)^{p'} dV_g(y)
  \right)^{\frac{1}{p'}}.
$$
Again by and Lemma \ref{l:lodlp} it follows that for $p > 2$ (and
hence for $p' < 2$) it is $JJJ \leq C r^{\frac 6p - 4}$. If we
choose $p \in (2, 3)$, then ${\frac 6p - 4} > - 2$, which implies
$JJJ < o_r(1) r^2$, and hence also
\begin{equation}\label{eq:v2l}
  \int_{A_{r,l}} v_{2,l} dV_g = o_r(1) r^2.
\end{equation}
Then, choosing first $R$ sufficiently large and then $l$
sufficiently large, \eqref{eq:dwlv1v2}, \eqref{eq:estv1l} and
\eqref{eq:v2l} conclude the proof.
\end{pf}

\subsection{Radial behavior}\label{ss:rb}

\

\noindent The next step consists in studying the dependence on $r$
of the function $\ov{u}_{r,l}$ defined in Lemma \ref{l:trueintw}.
It is well known that in geodesic coordinates the metric
coefficients $g_{ij}$ have the expression
\begin{equation}\label{eq:gij}
    g_{ij}(x) = \d_{ij} - \frac 13 R_{ikjl} x^k x^l + O(|x|^3),
\end{equation}
where $R_{ikjl}$ are the components of the curvature tensor, see
for example \cite{lp}, and the volume element satisfies
\begin{equation}\label{eq:dVg}
    dV_{g} = \sqrt{det g} \; dV_{\R^4} = ( 1 + O(|x|^2)) dV_{\R^4}
    \qquad \hbox{ with } \n \sqrt{det g} = O(|x|) \hbox{ and }
    \n^2 \sqrt{det g} = O(1).
\end{equation}
Using the exponential map at $x_l$, we can use coordinates $r,
\th$ in a neighborhood of $x_l$, where $r = |x| > 0$ and $\th \in
S^3$. In these coordinates the volume element $dV_g$ and the
surface element $d \s_g$ take the form
$$
  dV_g = r^3 \tilde{f}(r,\th) dr d\th; \qquad \qquad d\s_g =
  \tilde{f}(r,\th) d\th,
$$
where $\tilde{f}$ is a smooth bounded function on $\{r > 0\}$.
Using these coordinates, we consider a regular function $h$. Then,
letting $A_{\tilde{r}} = B_{2 \tilde{r}} (x_l) \setminus
B_{\tilde{r}} (x_l)$, one has
$$
  \int_{A_{\tilde{r}}} h dV_g = \int_{\tilde{r}}^{2 \tilde{r}} r^3 dr
  \int_{S^3} h(r,\th) f(r,\th) d \th; \qquad \frac{\partial h}{\partial
  \nu}(r,\th) = \frac{\partial h}{\partial r}(r,\th),
$$
where $\nu$ denotes the exterior unit normal to $\partial
B_{\tilde{r}}(x_l)$.

We also use the coordinates $z, \th$, where $z = \log r$. In these
new coordinates we obtain
$$
  dV_g = e^{4z} f(z,\th) dz d\th; \qquad \qquad d\s_g = e^{3z} f(z,\th) d\th,
$$
where $f(z,\th) = \tilde{f}(e^z,\th)$, and
$$
  \int_{A_{\tilde{r}}} h \; dV_g = \int_{s}^{s + \b} dz \int_{S^3} h(z,\th)
  f(z,\th) e^{4 z} d \th; \qquad \frac{\partial h}{\partial \nu}(z,\th)
  = e^{-z} \frac{\partial h}{\partial z}(z,\th).
$$
Here we have set $\b = \log 2$ and $s = \log \tilde{r}$. From
\eqref{eq:dVg} we also find
\begin{equation}\label{eq:dvgf}
    f(z,\th) = 1 + O(e^{2z}); \qquad \qquad \frac{\partial f}{\partial
    z}(z,\th) = O(e^{2z}); \qquad \qquad \frac{\partial^2 f}{\partial
    z^2}(z,\th) = O(e^{2z}).
\end{equation}
Now we can write
\begin{eqnarray}\label{eq:1stder}
    \frac{\partial}{\partial s} \int_{A_{\tilde{r}}} h \; dV_g & = & \left. \int_{S^3}
    h(z,\th) e^{4z} f(z,\th) d \th \right|_{z=s}^{s+\b} = \int_s^{s+\b}
    \int_{S^3} \frac{\partial}{\partial z} \left( h(z,\th) e^{4z} f(z,\th)
     \right) d \th d z \nonumber \\ & = & \int_s^{s+\b} \int_{S^3}
     \frac{\partial h}{\partial z} e^{4z} f(z,\th) d \th d z +
     \int_s^{s+\b} \int_{S^3} h(z,\th) \left( 4 f(z,\th) e^{4z} + e^{4z}
     \frac{\partial f}{\partial z} (z,\th) \right) d \th dz.
\end{eqnarray}
Taking a second derivative with respect to $s$, from the above
formulas we obtain
\begin{eqnarray*}
    \frac{\partial^2}{\partial s^2} \int_{A_{\tilde{r}}} h \; dV_g & = & \left.
    \int_{S^3} \frac{\partial h}{\partial z}(z,\th) e^{4z} h(z,\th) d \th
    \right|_{z=s}^{s+\b} + 4 \frac{\partial}{\partial s} \int_{A_{\tilde{r}}}
    h \; dV_g \\ & + & \frac{\partial}{\partial s} \left( \int_s^{s+\b} \int_{S^3}
    h(z,\th) e^{4z} \frac{\partial f}{\partial z} (z,\th) d\th dz\right)
    \\ & = & \int_{\partial A_{\tilde{r}}} e^{2 z} \frac{\partial h}{\partial \nu} d
    \s_g + 4 \frac{\partial}{\partial s} \int_{A_{\tilde{r}}}
    h \; dV_g + \frac{\partial}{\partial s} \left( \int_s^{s+\b} \int_{S^3}
    h(z,\th) e^{4z} \frac{\partial f}{\partial z} (z,\th) d\th dz\right).
\end{eqnarray*}
Using the coordinates $(r,\th)$ and integrating by parts we derive
\begin{eqnarray*}
  \int_{\partial A_{\tilde{r}}} e^{2 z} \frac{\partial h}{\partial \nu} d
    \s_g & = & \int_{\partial A_{\tilde{r}}} r^2 \frac{\partial h}{\partial \nu} d
    \s = \int_{A_{\tilde{r}}} r^2 \D h \; dV_g - \int_{A_{\tilde{r}}} h \D r^2 dV_g
    + \int_{\partial A_{\tilde{r}}} h \frac{\partial r^2}{\partial \nu} d \s_g \\
    & = & \int_{A_{\tilde{r}}} r^2 \D h \;
    dV_g - 8 \int_{A_{\tilde{r}}} h \; dV_g + 2 \int_{\partial A_{\tilde{r}}}
    h e^{4 z} d \s_g + \int_{A_{\tilde{r}}} (\D r^2 - 8) h \; dV_g.
\end{eqnarray*}
By the last two formulas we finally get the following equation
\begin{eqnarray}\label{eq:diffe}
    \frac{\partial^2}{\partial s^2} \int_{A_{\tilde{r}}} h \; dV_g & = &
    6 \frac{\partial}{\partial s} \int_{A_{\tilde{r}}} h \; dV_g - 8
    \int_{A_{\tilde{r}}} h \; dV_g + \int_{A_{\tilde{r}}} r^2 \D h \; dV_g \\
    & + & \int_{A_{\tilde{r}}} (\D r^2 - 8) h \; dV_g
    + \frac{\partial}{\partial s} \left( \int_s^{s+\b} \int_{S^3}
    \nonumber
    h(z,\th) e^{4z} \frac{\partial f}{\partial z} (z,\th) d\th dz\right).
\end{eqnarray}
Next we want to apply \eqref{eq:diffe} to the case of $h = u_l$,
and derive a differential equation involving the average
$\ov{u}_{r,l}$ of $u_l$ on the annuli $A_{r,l}$.

\begin{lem}\label{l:system}
Suppose that $(x_l)_l \subseteq M$, $(s_l)_l \subseteq \R_+$,
$i(M) \geq s_l > 0$, and that $(u_l)_l$ satisfies \eqref{eq:pl}
and \eqref{eq:noul}. Then, for every $l$ and every $r < s_l$ we
let
$$
  W_l(z) = \frac{1}{Vol(A_{r,l})} \int_{A_{r,l}} u_l \; dV_g;
 \qquad \qquad z = \log r,
$$
where $A_{r,l}$ is defined in \eqref{eq:}. Then the functions
$W_l(z)$ solve the following equation
\begin{equation}\label{eq:system}
    W''_l(z) + 2 (1 + O(e^{2z})) W'_l(z) = \frac{\int_{A_{r,l}} r^2
    \D_g u_l d V_g}{Vol(A_{r,l})} + O(e^{2z}), \qquad \hbox{ for } z \in
    (\log (a_l r_l), \log s_l).
\end{equation}
\end{lem}

\noindent We notice that the function $W_l(z)$ coincide with
$\ov{u}_{r,l}$ up to the change of variables $r \mapsto z = \log
r$.

\

\begin{pf}
We first let
$$
  \tilde{W}_l(z) = \int_{A_{r,l}} u_l \, dV_g; \quad \qquad Y_l(r) =
  \int_{A_{r,l}} dV_g, \qquad \qquad \quad z = \log r.
$$
We have clearly
$$
  W'_l(z) = \left( \frac{\tilde{W}_l(z)}{Y_l(z)} \right)' =
  \frac{\tilde{W}'_l(z) Y_l(z) - Y'_l(z)
  \tilde{W}_l(z)}{Y_l^2(z)},
$$
and
$$
  W''_l(z) = \frac{Y_l^2(z)
  \left[ \tilde{W}''_l(z) Y_l(z) - Y''_l(z) \tilde{W}_l(z) \right] -
  2 Y_l(z) Y'_l(z) \left[ \tilde{W}'_l(z) Y_l(z) - Y'_l(z) \tilde{W}_l(z)
  \right]}{Y_l^4(z)}.
$$
Using the last two formulas and \eqref{eq:diffe} with
$A_{\tilde{r}} = A_{r,l}$ and $h = u_l$, after some calculation
(which also uses \eqref{eq:1stder} with $h$ replaced by
$\frac{h}{f} \frac{\partial f}{\partial z}$) we obtain
\begin{eqnarray*}
    W''_l(z) & = & 6 W'_l(z) - 2 \frac{Y'_l(z)}{Y_l(z)} W'_l(z) +
    \frac{\int_{A_{r,l}} r^2 \D_g u_l d V_g}{Y_l(z)} \\ & + & \left[
    \int (\D_g r^2 - 8) u_l + \int \frac{\partial}{\partial z} \left(
    \frac{u_l \frac{\partial f}{\partial z}}{f} \right) e^{4z} f +
    \int u_l \frac{\frac{\partial f}{\partial z}}{f} \left( 4 f e^{4z}
    + e^{4z} \frac{\partial f}{\partial z} \right) \right]
    \frac{\int e^{4z} f}{Y_l(z)^2} \\ & - & \left[
    \int (\D_g r^2 - 8) + \int \frac{\partial}{\partial z} \left(
    \frac{\frac{\partial f}{\partial z}}{f} \right) e^{4z} f +
    \int \frac{\frac{\partial f}{\partial z}}{f} \left( 4 f e^{4z}
    + e^{4z} \frac{\partial f}{\partial z} \right) \right]
    \frac{\int u_l e^{4z} f}{Y_l(z)^2}.
\end{eqnarray*}
We notice that, adding and subtracting the average of
$\ov{u}_{r,l}$ to $u_l$, some cancellation occurs. Moreover, from
\eqref{eq:dvgf} and \eqref{eq:1stder} we get
$$
  \frac{Y'_l(z)}{Y_l(z)} = \frac{\int \left( 4 e^{4z} f + e^{4z}
  \frac{\partial f}{\partial z} \right)}{Y_l(z)} = 4 + O(e^{2z}).
$$
Therefore, using these remarks we obtain
\begin{eqnarray*}
    W''_l(z) & = & - 2 (1 + O(e^{2z})) W'_l(z) +
    \frac{\int_{A_{r,l}} r^2 \D_g u_l d V_g}{Y_l(z)} \\ & + & \left[
    \int (\D_g r^2 - 8) (u_l - \ov{u}_{r,l}) + \int \frac{\partial}{\partial z}
    \left( \frac{(u_l - \ov{u}_{r,l}) \frac{\partial f}{\partial z}}{f} \right)
    e^{4z} f \right. \\ & + & \left. \int (u_l - \ov{u}_{r,l})
    \frac{\frac{\partial f}{\partial z}}{f} \left( 4 f e^{4z}
    + e^{4z} \frac{\partial f}{\partial z} \right) \right]
    \frac{\int e^{4z} f}{Y_l(z)^2} \\ & - & \left[
    \int (\D_g r^2 - 8) + \int \frac{\partial}{\partial z} \left(
    \frac{\frac{\partial f}{\partial z}}{f} \right) e^{4z} f +
    \int \frac{\frac{\partial f}{\partial z}}{f} \left( 4 f e^{4z}
    + e^{4z} \frac{\partial f}{\partial z} \right) \right]
    \frac{\int (u_l - \ov{u}_{r,l}) e^{4z} f}{Y_l(z)^2}.
\end{eqnarray*}
We next estimate the terms in the last three lines of this
expression. We begin by noticing that $(\D r^2 - 8) = O(r^2)$,
which can be deduced from elementary computations in local
coordinates. This and the Poincar\'e inequality imply
$$
  \left| \int (\D_g r^2 - 8) (u_l - \ov{u}_{r,l}) d V_g \right| \leq C
  e^{3z} \int_{A_{r,l}} |\n u_l| d V_g, \qquad \quad z = \log r.
$$
From Lemma \ref{l:lodlp} then one finds
$$
 \left| \int (\D_g r^2 - 8) (u_l - \ov{u}_{r,l}) d V_g \right| \leq C
  e^{6z}.
$$
Similarly, using \eqref{eq:dvgf} and also the fact that
$\frac{\partial u_l}{\partial z} = \frac{\partial u_l}{\partial r}
\frac{\partial r}{\partial z} = O(e^{z} |\n u_l|)$, we obtain
\begin{eqnarray*}
    \left| \int \frac{\partial}{\partial z} \left( \frac{(u_l - \ov{u}_{r,l})
    \frac{\partial f}{\partial z}}{f} \right) e^{4z} f \right| & \leq &
    \int_{A_{r,l}} O(e^{2z}) |u_l - \ov{u}_{r,l}| d V_g + \int_{A_{r,l}}
    O(e^{3z}) |\n u_l| d V_g \\ & \leq & C e^{6z}.
\end{eqnarray*}
Reasoning in the same way for the remaining terms we finally
deduce
$$
  W''_l(z) + 2 (1 + O(e^{2z})) W'_l(z) = \frac{\int_{A_{r,l}} r^2 \D_g h d
  V_g}{Y_l(z)} + O(e^{2z}).
$$
Then the last four estimates imply the first equation in
\eqref{eq:system}.
\end{pf}

\begin{rem}\label{r:bdW'l}
Using \eqref{eq:1stder} with $A_{\tilde{r}} = A_{r,l}$, and with
$h = u_l$ (or with $h = 1$ to compute $Y'_l$), we obtain
\begin{eqnarray*}
    W'_l(z) = \frac{\left[ \int_{A_{r,l}} u_l \left( 4 f e^{4z} + e^{4z}
    \frac{\partial f}{\partial z} e^{4z} f \right) + \int_{A_{r,l}} \frac{\partial u_l}{\partial z}
    \right] \int_{A_{r,l}} f e^{4z} - \left[ \int_{A_{r,l}} 4 e^{4z} f + e^{4z} \frac{\partial
    f}{\partial z} \right] \int_{A_{r,l}} u_l f e^{4z}}{\left( \int_{A_{r,l}} f e^{4z}
    \right)^2}.
\end{eqnarray*}
If we denote again by $\ov{u}_{r,l}$ the average of $u_l$ in the
annulus $A_{r,l}$, adding and subtracting this average from $u_l$
in the last formula we get some cancellations and we are left with
\begin{eqnarray*}
    W'_l(z) & = & \frac{\left[ \int_{A_{r,l}} (u_l - \ov{u}_{r,l}) \left(
  e^{4z} \frac{\partial f}{\partial z} \right) +
  \int_{A_{r,l}} \frac{\partial u_l}{\partial z} e^{4z} f \right] \int_{A_{r,l}}
  f e^{4z}}{\left( \int_{A_{r,l}} f e^{4z} \right)^2} -
    \frac{\left[ \int_{A_{r,l}} e^{4z} \frac{\partial
    f}{\partial z} \right] \int_{A_{r,l}} (u_l - \ov{u}_{r,l})
    f e^{4z}}{\left( \int_{A_{r,l}} f e^{4z} \right)^2}.
\end{eqnarray*}

As a byproduct of this formula and the Poincar\'e inequality we
deduce
$$
  |W'_l(z)| \leq C \frac{\int_{A_{r,l}} |u_l - \ov{u}_{r,l}| d
  V_g}{Y_l(z)} + C r \frac{\int_{A_{r,l}} |\n u_l| d V_g}{Y_l(z)}
  \leq C r \frac{\int_{A_{r,l}} |\n u_l| d V_g}{Y_l(z)}.
$$
Then, applying Lemma \ref{l:lodlp}, we find
\begin{equation}\label{eq:bdW'l}
  |W'_l(z)| \leq C.
\end{equation}
\end{rem}

\

\noindent In the next lemma we study the solutions of
\eqref{eq:system} in the case of a simple blow-up. When $x_0 = 0$
and $\l = 1$, the function $\hat{u}_\infty$, see \eqref{eq:hwinf},
is of the form
$$
  \hat{u}_\infty(x) = \log \left( \frac{2}{1 + |x|^2} \right) +
  \frac 14 \log \frac{3}{k_0}.
$$
From straightforward computations one finds
\begin{eqnarray*}
    \int_{A_r} \hat{u}_\infty dV_{\R^4} & = & 2 \pi^2 \left[  \frac{15}{4}
    r^4 \log 2 + 4r^4 \log \left( \frac{1}{1+4r^2} \right) + \frac{15}{8}
    r^4 - \frac 34 r^2 + \frac 14 \log(1+4r^2) \right. \\ & - & \left.
    \frac 14 r^4 \log \left( \frac{1}{1+r^2} \right) - \frac 14 \log
    (1+r^2) \right].
\end{eqnarray*}
Scaling back to $u_l$, using \eqref{eq:rlslxl} and some elementary
estimates one deduces (for $t > 0$ large and fixed)
\begin{equation}\label{eq:icAl}
    W_l(\log r_l + t) = - 2 t + \ov{C} - \log r_l + O(e^{-2t}) +
  o_l(1); \qquad W'_l(\log r_l + t) = - 2 + O(e^{-2t}) + o_l(1),
\end{equation}
where $\ov{C}$ is some explicit positive constant.

Now we prove some upper bounds for the function $W_l$. Notice from
\eqref{eq:icAl} that $W_l$ at $z = \log r_l + t$ ($t$ large and
fixed) has slope close to $-2$. Given $\g \in (1,2)$, we consider
an affine function $h^\g_{t,l}$ which coincides with $W_l$ for $z
\sim \log r_l$ and which has slope $- \g > -2$. The next lemma
asserts that indeed $W_l(z) < h^\g_{t,l}(z)$ until $z$ gets close
to $\log s_l$. This is helpful to get integral estimates on $e^{4
u_l}$, which is done at the end of the section.

\begin{lem}\label{l:ODE}
Suppose $(x_l)_l$, $(r_l)_l$, $(s_l)_l$ are a simple blow-up for
$(u_l)_l$, and let $(W_l)_l$ be given by Lemma \ref{l:system}.
Given $\g \in (1,2)$ and $t > 0$, consider the following functions
$$
  h^\g_{t,l}(z) = - \g (z - \log r_l - t) + W_l(\log r_l + t).
$$
Then there exist $t_l \to + \infty$ arbitrarily slowly and $C_\g
> 0$ such that for $l$ large
$$
  W_l(z) \leq h_{t_l,l}^\g(z); \qquad \qquad z \in \left[
  \log r_l + t_l, \log s_l - C_\g \right].
$$
\end{lem}

\begin{pf}
Recall that $(W_l)_l$ are solutions of \eqref{eq:system}
satisfying the {\em initial} conditions \eqref{eq:icAl} for any
large and fixed $t$. If $t_l \to + \infty$ sufficiently slowly, we
can also replace $t$ by $t_l$ in \eqref{eq:icAl}, namely we can
also assume that
\begin{equation}\label{eq:unif}
    W_l(\log r_l + t_l) = - 2 t_l + \ov{C} - \log r_l + o_l(1);
    \qquad \qquad W'_l(\log r_l + t_l) = - 2 + o_l(1).
\end{equation}

Suppose by contradiction that there exist $\ov{s}_l \in [\log r_l,
\log s_l]$, with $\log s_l - \ov{s}_l \to + \infty$ such that
$W_l$ intersects $h_{t_l,l}^\g$ for the first time. We notice
that, by the asymptotics in \eqref{eq:icAl}, it must also be
$\ov{s}_l - \log r_l - t_l \to + \infty$ if $t_l \to + \infty$
sufficiently slowly. Then we have
$$
  W_l(\ov{s}_l) = h^{\g}_{t_l,l}(\ov{s}_l);
  \qquad \qquad \qquad W'_l(\ov{s}_l) \geq - \g.
$$
We now choose a sequence of real numbers $(H_l)_l$ by means of the
following condition
$$
  H_l = \sup \left\{ H \in \R \; : \; h_{t_l,l}^{\frac{\g+2}{2}} +
  H < W_l \quad \hbox{ in } [\log r_l + t_l, \ov{s}_l] \right\}.
$$
By \eqref{eq:icAl} it must be $H_l \to - \infty$ as $l \to +
\infty$ (provided $t_l \to + \infty$ sufficiently slowly), and
there exist $\tilde{s}_l$ such that
\begin{equation}\label{eq:tsl1}
  W_l(\tilde{s}_l) = h_{t_l,l}^{\frac{\g+2}{2}}(\tilde{s}_l) +
  H_l; \qquad \qquad W'_l(\tilde{s}_l) = - \frac{\g+2}{2}; \qquad
  \qquad W''_l(\tilde{s}_l) \geq 0.
\end{equation}
Moreover, by \eqref{eq:bdW'l} and \eqref{eq:icAl}, $\tilde{s}_l$
satisfies
\begin{equation}\label{eq:tsl2}
    |\ov{s}_l - \tilde{s}_l| \to + \infty \hbox{ as } l \to +
    \infty; \qquad \qquad |\tilde{s}_l - \log r_l - t_l| \to +
    \infty \hbox{ as } l \to + \infty.
\end{equation}

Next we claim that, for $C > 0$ sufficiently large, the following
property holds
\begin{equation}\label{eq:cla1}
  \int_{B_{\frac{e^{\ov{s}_l}}{C}}(x_l) \setminus B_{e^{t_l} r_l}(x_l)}
  e^{4 u_l} d V_g \to 0 \qquad \hbox{ as } l \to + \infty.
\end{equation}
In order to prove this claim, let us recall that by our choice of
$\ov{s}_l$, it is $W_l(z) \leq h^\g_{t_l,l}(z)$ for every $z \in
[\log r_l + t_l, \ov{s}_l]$. Dividing the region
$B_{\frac{e^{\ov{s}_l}}{C}}(x_l) \setminus B_{e^{t_l} r_l}(x_l)$
into concentric annuli $A'_{r,l}$ (see \eqref{eq:2}) of suitable
radii, we get
$$
  \int_{B_{\frac{e^{\ov{s}_l}}{C}}(x_l) \setminus B_{e^{t_l} r_l}(x_l)}
  e^{4 u_l} dV_g \leq
  \sum_{j=0}^{j_l} \int_{A'_{\hat{r}_{l,j},l}} e^{4 u_l} dV_g,
$$
where
$$
  \hat{r}_{l,j} = \frac 45 e^{t_l} r_l \left( \frac 75 \right)^j ; \qquad
  \qquad \qquad \left( \frac 75 \right)^{j_l} \in \left( \frac 54
  \frac{e^{\ov{s}_l}}{C e^{t_l} r_l}, \frac 52 \frac{e^{\ov{s}_l}}{C e^{t_l} r_l}
  \right).
$$
Given $\g \in (1,2)$, from Lemma \ref{l:trueintw} it follows that
\begin{eqnarray*}
    \int_{A'_{\hat{r}_{l,j},l}} e^{4 u_l} dV_g \leq C
    |A_{\hat{r}_{l,j},l}| e^{4 \ov{u}_{l,\hat{r}_{l,j}}} \leq C
    \hat{r}_{l,j}^4 e^{4 W_l(\log \hat{r}_{l,j})} \leq C
    \hat{r}_{l,j}^4 e^{4 h_{t_l,l}^\g(\log \hat{r}_{l,j})}; \qquad
    j = 1, \dots, j_l.
\end{eqnarray*}
From the expression of $h^\g_{t_l,l}$ and \eqref{eq:unif} we
deduce
\begin{eqnarray*}
    \hat{r}_{l,j}^4 e^{4 h_{t_l,l}^\g(\log \hat{r}_{l,j})} & \leq &
    C \hat{r}_{l,j}^4 \exp \left[ 4 \left( - \g \left( \log \hat{r}_{l,j}
    - \log r_l - t_l \right) - 2 t_l + \ov{C} - \log r_l + o_l(1) \right)
    \right] \\ & = & C \hat{r}_{l,j}^4 \exp \left[ - 4 \g \log \hat{r}_{l,j}
    + 4 (\g - 1) \log r_l + 4 (\g - 2) t_l + \ov{C} + o_l(1) \right]
    \\ & \leq & C \left( \frac{r_l}{\hat{r}_{l,j}}
    \right)^{4(\g-1)} e^{4(\g-2)t_l} = C \left( \frac{5}{4 e^{t_l}}
    \right)^{4(\g-1)} e^{4(\g-2)t_l} \left( \frac 57 \right)^{4(\g-1)j}.
\end{eqnarray*}
Hence it follows that
$$
  \int_{B_{\frac{e^{\ov{s}_l}}{C}}(x_l) \setminus B_{e^{t_l} r_l}(x_l)} e^{4 u_l}
  dV_g \leq C \left( \frac{5}{4 e^{t_l}} \right)^{4(\g-1)} e^{4(\g-2)t_l}
  \sum_{j=0}^{\infty} \left( \frac 57 \right)^{4(\g-1)j} \to 0,
$$
since $\g \in (1,2)$ and since $t_l \to + \infty$. This proves
\eqref{eq:cla1}.

We can now apply Lemma \ref{l:harlap} with $\Sig_l = e^{t_l} r_l$,
$S_l = \frac{e^{\ov{s}_l}}{C}$, and $\log r = \tilde{s}_l$. Also,
by \eqref{eq:tsl2} and \eqref{eq:cla1}, we can choose $\e = \e_l
\to 0$ and $R = R_l \to + \infty$ sufficiently slowly. Therefore,
from Lemma \ref{l:harlap} and Proposition \ref{p:hwlb} (see in
particular \eqref{eq:intalrlu}) we deduce that
$$
  \int_{A_{e^{\tilde{s}_l},l}} |x-x_l|^2 (- \D u_l(x)) dV_g(x) = \left(
  \frac{15}{8} \int_{B_{\frac{e^{\tilde{s}_l}}{R_l}}(x_l)} f_l dV_g + o_l(1)
  \right) e^{4 \tilde{s}_l} \geq \left( 30 \pi^2 + o_l(1)
  \right) e^{4 \tilde{s}_l}.
$$
On the other hand, from \eqref{eq:system} and the last two
conditions in \eqref{eq:tsl1} we find
\begin{eqnarray*}
    \int_{A_{e^{\tilde{s}_l},l}} |x-x_l|^2 (- \D u_l(x)) dV_g(x) & =
    & \left[ - W''_l(\tilde{s}_l) - 2 (1 + O(e^{2 \tilde{s}_l}))
    W'_l(\tilde{s}_l) + O(e^{2 \tilde{s}_l}) \right]
    Y_l(\tilde{s}_l) \\ & \leq & \left[ \g + 2 + o_l(1) \right]
    \left( \frac{15 \pi^2}{2} + o_l(1) \right) e^{4
    \tilde{s}_l}.
\end{eqnarray*}
Since $\g < 2$, from the last two inequalities we get a
contradiction. This concludes the proof of the Lemma.
\end{pf}

\

\noindent We are finally in position to prove Proposition
\ref{p:sbu}.

\

\begin{pfn} {\sc of Proposition \ref{p:sbu}.}
It is sufficient to apply Lemma \ref{l:ODE} and to reason as for
the proof of \eqref{eq:cla1}. In fact, in this way we get
$$
  \int_{B_{\frac{e^{{s}_l}}{C}}(x_l) \setminus B_{e^{t_l}
  r_l}(x_l)} e^{4 u_l} d V_g \to 0 \qquad \hbox{ as } l \to + \infty.
$$
Moreover, choosing $b_l = e^{t_l}$ in \eqref{eq:intalrlu} and $t_l
\to + \infty$ sufficiently slowly, we also get
$$
    \int_{B_{b_l r_l}(x_l)} e^{4 u_l} dV_g \to \frac{8 \pi^2}{k_0}
    \qquad \hbox{ as } l \to + \infty.
$$
The last two formulas yield the conclusion.
\end{pfn}

\section{Proof of Theorem \ref{th:bd}}\label{s:fa}

We prove first the theorem under the assumption \eqref{eq:kp3},
and we postpone the remaining cases to a second subsection.

\subsection{Proof under the assumption \eqref{eq:kp3}}\label{ss:kp3}

In this subsection we show how a general blow-up phenomenon can be
essentially reduced to the case of finitely-many simple blow-ups.
We divide the proof into three steps, and we always assume that
$(u_l)_l$ is a sequence satisfying \eqref{eq:noul} and
\eqref{eq:unbdul}. We recall that the integer $k$ is defined by
the condition $k_0 \in (8 k \pi^2, 8 (k+1) \pi^2)$.

\

\noindent {\bf Step 1.} There exist an integer $j \leq k$,
sequences $(x_{1,l})_l, \dots (x_{j,l})_l \subseteq M$ and radii
$(r_{1,l})_l, \dots, (r_{j,l})_l$, $(\tilde{r}_{1,l})_l, \dots,
(\tilde{r}_{j,l})_l \to 0$ satisfying the properties (for some $\a
\in (0,1)$)
\begin{equation}\label{eq:step11}
   \frac{\tilde{r}_{i,l}}{r_{i,l}} \to + \infty \hbox{ (slowly)}
   \hbox{ as } l \to + \infty; \qquad \qquad B_{\tilde{r}_{i,l}} \cap
   B_{\tilde{r}_{h,l}} = \emptyset \hbox{ for } i \neq h;
\end{equation}
\begin{equation}\label{eq:step12}
  \forall R > 0 \; \; \hat{u}_{l,i}
    \to \log \frac{2}{1 + |x|^2} - \frac{1}{4} \log \left( \frac 13 k_0
    \right) \hbox{ in } H^4(B_R^{\R^4}) \cap C^\a(B_R^{\R^4})
    \hbox{ as } l \to + \infty.
\end{equation}
\begin{equation}\label{eq:step13}
   \forall \rho > 0 \exists C_\rho > 0 \hbox{ s.t. if } \int_{B_s(y)} e^{4 u_l}
    d V_g \geq \rho \hbox{ with } B_s(y) \subseteq
    M \setminus \cup_{i=1}^j B_{\tilde{r}_{i,l}}(x_{i,l}), \hbox{ then } s \geq
    C_\rho^{-1} d_l(y),
\end{equation}
where $d_l(y) = \min_{i=1, \dots, j} |y - x_{i,l}|$. Here
$\hat{u}_{l,j}$ denotes the function obtained using the procedure
in Section \ref{s:bub}, but scaling around the point $x_{i,l}$
with a factor $r_{i,l}$.

\

\noindent In order to prove Step 1, we consider a small number
$\rho > 0$, say $\rho \in \left( 0, \frac{\pi^2}{k_0} \right)$,
and we define sequences $(x_{1,l})_l \subseteq M$, $(r_{1,l})_l
\subseteq \R_+$ satisfying
$$
  \int_{B_{r_{1,l}}(x_{1,l})} e^{4 u_l} dV_g = \max_{x \in M}
  \int_{B_{r_{1,l}}(x)} e^{4 u_l} dV_g = \rho.
$$
If \eqref{eq:unbdul} holds, it must be $r_{1,l} \to 0$ as $l \to +
\infty$. In fact, if it were $r_{1,l} \geq C^{-1}$, we could apply
Proposition \ref{p:cc} to get uniform $L^p$ bounds on $e^{4(u_l -
\ov{u}_l)}$ for some $p > 1$. This fact and the Jensen inequality
would yield
$$
   1 = e^{4 \ov{u}_l} \int_M e^{4(u_l - \ov{u}_l)} dV_g \leq C
   e^{4 \ov{u}_l}; \qquad \quad \ov{u}_l \leq C,
$$
and hence uniform bounds on $e^{4 u_l}$ in $L^p(M)$. This would
imply, by elliptic regularity results, uniform bounds in $H^2(M)$
on $(u_l)_l$, which is a contradiction to our assumptions.

Then, if $\frac{\tilde{r}_{1,l}}{r_{1,l}}$ tends to infinity
sufficiently slowly, $(r_{1,l})_l$ and $(\tilde{r}_{1,l})_l$
satisfy \eqref{eq:trlxl}, so Proposition \ref{p:hwlb} applies
yielding the existence of a bubble, giving \eqref{eq:step12} for
$i = 1$ and
$$
  \int_{B_{\tilde{r}_{1,l}}(x_{1,l})} e^{4 u_l} dV_g = \frac{8
  \pi^2}{k_0} + o_l(1).
$$

\noindent If \eqref{eq:step13} holds for $j = 1$, Step 1 is
proved.

\

\noindent If \eqref{eq:step13} does not hold, there exists $\rho_1
> 0$, which can be assumed belonging to $\left( 0,
\frac{\pi^2}{k_0}\right)$, and there exist sequences $(y_l)_l
\subseteq M$, $\ov{r}_l \subseteq \R_+$ such that
\begin{equation}\label{eq:rho1}
    \int_{B_{\ov{r}_l}(y_l)} e^{4 u_l} dV_g \geq \rho_1; \qquad
    B_{\ov{r}_l}(y_l) \subseteq M \setminus
    B_{\tilde{r}_{1,l}}(x_{1,l}); \qquad \frac{\ov{r}_l}{|y_l -
    x_{1,l}|} \to 0 \quad \hbox{ as } l \to + \infty.
\end{equation}
Now we define $r_{2,l}$ and $x_{2,l}$ such that
$$
  \int_{B_{r_{2,l}}(x_{2,l})} e^{4 u_l} dV_g = \max_{B_{r_{2,l}}(y)
  \subseteq M \setminus B_{\tilde{r}_{1,l}}(x_{1,l})} \int_{B_{r_{2,l}}(y)}
  e^{4 u_l} dV_g = \rho_1.
$$
By Proposition \ref{p:hwlb} it is easy to see that if
$\frac{\tilde{r}_{1,l}}{r_{1,l}} \to + \infty$ sufficiently
slowly, then we have
\begin{equation}\label{eq:rhrh2}
    \frac{\tilde{r}_{1,l}}{|x_{1,l} - x_{2,l}|} \to 0; \qquad
   \qquad \frac{r_{2,l}}{|x_{1,l} - x_{2,l}|} \to 0 \qquad \qquad
   \qquad \hbox{ as } l \to + \infty,
\end{equation}
which in particular implies $r_{2,l} \to 0$ as $l \to + \infty$.
Therefore, by the last formula we can find $\hat{r}_{2,l}
\subseteq \R_+$ such that
$$
   \int_{B_{r_{2,l}}(y)} e^{4 u_l} d V_g \leq \rho_1 \hbox{ for every }
   y \in B_{\hat{r}_{2,l}}(x_{2,l}); \qquad \qquad
   \frac{\hat{r}_{2,l}}{|x_{1,l} - x_{2,l}|} \to 0
   \qquad \hbox{ as } l \to + \infty.
$$
Then Proposition \ref{p:hwlb} applies yielding the existence of a
second bubble.

Continuing in this way, we see immediately that $j$ cannot exceed
$k$, since every bubble contributes an amount of $\frac{8
\pi^2}{k_0}$ to the volume and since we are assuming
\eqref{eq:noul}. This concludes the proof of Step 1.

\

\noindent {\bf Step 2.} If in Step $1$ it is $j = 1$, then there
holds
\begin{equation}\label{eq:star4}
    \int_M e^{4 u_l} dV_g = \frac{8 \pi^2}{k_0} + o_l(1).
\end{equation}

\

\noindent In this case, if we choose $s_l = \frac 12 i(M)$ for
every $l$, where $i(M)$ is the injectivity radius of $M$, then by
\eqref{eq:step13} $(x_{1,l})_l$, $(r_{1,l})_l$, $(s_l)_l$ are a
simple blow-up for $u_l$. Therefore Proposition \ref{p:sbu}
applies and, since $(s_l)_l$ is uniformly bounded from below,
there exists $C > 0$ such that for $l$ large
\begin{equation}\label{eq:michi}
    \int_{B_{C^{-1}}(x_{1,l})} e^{4 u_l} = \frac{8 \pi^2}{k_0}
  + o_l(1).
\end{equation}
We prove first the following property
\begin{equation}\label{eq:clst2}
    \ov{u}_l \to - \infty \qquad \hbox{ as } l \to + \infty.
\end{equation}
In fact, using the Green's representation formula, for a.e. $x \in
M$ we obtain
$$
  u_l(x) = \ov{u}_l + \int_M G(x,y) \left( 2 k_l e^{4 u_l}(y) - 2 Q_l
  \right) d V_g(y) \geq \ov{u}_l - C + \int_M G(x,y) 2 k_l e^{4 u_l}(y)
  d V_g(y).
$$
By \eqref{eq:step12} and \eqref{eq:intehwi}, given any small
$\tilde{\e} > 0$, there exists $R_{\tilde{\e}}$ such that, for $l$
sufficiently large
$$
  \int_{B_{R_{\tilde{\e}} r_{1,l}}(x_{1,l})} 2 k_l e^{2 u_l} \geq 16 \pi^2
  - 2 \pi^2 \tilde{\e}.
$$
Hence the last two formulas and \eqref{eq:G} imply
$$
  e^{4 u_l(x)} \geq C^{-1} e^{4 \ov{u}_l} \frac{1}{|x -
  x_{1,l}|^{8-\tilde{\e}}}; \qquad \qquad \hbox{ for } |x - x_{1,l}|
  \geq 2 R_{\tilde{\e}} r_{1,l},
$$
from which it follows that
\begin{equation}\label{eq:eqeq3}
    \int_M e^{4 u_l} d V_g \geq \int_{B_{i(M)}(x_{1,l}) \setminus
  B_{2 R_{\tilde{\e}} r_{1,l}}(x_{1,l})} e^{4 u_l} dV_g \geq
  C^{-1} e^{4 \ov{u}_l} \int_{2 R_{\tilde{\e}} r_{1,l}}^{i(M)}
  s^{\tilde{\e} - 5} ds \geq C^{-1} e^{4 \ov{u}_l}
  (R_{\tilde{\e}} r_{1,l})^{\tilde{\e} - 4}.
\end{equation}
If $\tilde{\e}$ is sufficiently small, the last factor tends to $+
\infty$ as $l \to + \infty$. Therefore \eqref{eq:clst2} follows
from \eqref{eq:noul}.

Now, by \eqref{eq:step13}, we can cover $M \setminus
B_{C^{-1}}(x_{1,l})$ with a finite number of balls $B_{r_i}(y_i)$,
$i = 1, \dots, \ell$ such that for every $i$ there holds
$\int_{B_{2r_i}(y_i)} e^{4 u_l} d V_g \leq \frac{\pi^2}{k_0}$.
Reasoning as in the proof of Proposition \ref{p:cc} one then finds
\begin{eqnarray*}
    \int_{M \setminus B_{C^{-1}}(x_{1,l})} e^{4 u_l} \leq C
    e^{4 \ov{u}_l} \sup_{y \in M, i = 1, \dots, \ell} \int_M \left(
    \frac{1}{|x - y|} \right)^{\frac{4
    \|e^{4 u_l}\|_{L^1(B_{2 r_i}(y_i))}}{8 \pi^2}} \leq C e^{4 \ov{u}_l} \to 0.
\end{eqnarray*}
Then \eqref{eq:michi} and the last formula conclude the proof of
Step 2.

\

\noindent {\bf Step 3.} If $j$ in Step $1$ is arbitrary, there
holds
\begin{equation}\label{eq:star44}
    \int_M e^{4 u_l} dV_g = \frac{8 \pi^2}{k_0} j + o_l(1).
\end{equation}

\

\noindent If $j > 1$ we reason as in \cite{ls}, and we analyze the
clustering of accumulation points. By relabelling the indices, we
can assume that
\begin{equation}\label{eq:x1lx2l}
    |x_{1,l} - x_{2,l}| = \inf_{i \neq h} |x_{i,l} - x_{h,l}| \to 0
    \hbox{ as } l \to + \infty.
\end{equation}
Of course, if $\inf_{i \neq h} |x_{i,l} - x_{h,l}| \not\to 0$,
then we could reason as in Step 2 a finite number of times.
Assuming \eqref{eq:x1lx2l}, we consider the set $X_{1,l} \subseteq
\{x_{1,l}, \dots, x_{h,l}\}$ of accumulation points for which the
distance from $x_{1,l}$ is comparable to $|x_{1,l} - x_{2,l}|$,
namely for which there exists $C > 0$ (independent of $l$) such
that
$$
  |x_{i,l} - x_{1,l}| \leq C |x_{1,l} - x_{2,l}|; \qquad i = 2,
  \dots, h = card(X_{1,l}).
$$
By our choices of the points $x_{1,l}, \dots, x_{h,l}$ and by
\eqref{eq:x1lx2l}, one easily checks that the three sequences
$(x_{i,l})_l, (r_{i,l})_l$ and $C^{-1} |x_{1,l} - x_{2,l}|$, $i =
1, \dots, h$, are a simple blow-up if $C$ is sufficiently large,
and Proposition \ref{p:sbu} applies yielding
\begin{equation}\label{eq:eded}
    \int_{B_{C^{-1} |x_{1,l} - x_{2,l}|}(x_{i,l})} e^{4 u_l}
    dV_g = \frac{8 \pi^2}{k_0} + o_l(1); \qquad i = 1, \dots, h.
\end{equation}

\

\noindent Our next claim is that there is no further concentration
in a neighborhood of $X_{1,l}$ of size comparable to $|x_{1,l} -
x_{2,l}|$. More precisely we have the following result.

\begin{lem}\label{l:step31}
In the above notation, for any large and fixed $C$ there holds
\begin{equation}\label{eq:star2}
  \int_{B_{C |x_{1,l} - x_{2,l}|}(x_{1,l})} e^{4 u_l} dV_g =
  \frac{8 \pi^2}{k_0} card(X_{1,l}) + o_l(1).
\end{equation}
\end{lem}

\begin{pf} In order to prove this claim we use a variant of the
argument in Step 2. First of all, for $\rho$ small and fixed, we
can cover the set $B_{C |x_{1,l} - x_{2,l}|}(x_{1,l}) \setminus
\cup_{i = 1, \dots, h} B_{C^{-1} |x_{1,l} - x_{2,l}|}(x_{i,l})$
with $\ell_l$ balls $B_{\rho_{n,l}}(y_{n,l})$, $n = 1, \dots,
\ell_l$, with the following properties
\begin{equation}\label{eq:rhoil}
    \ell_l \leq C; \qquad C^{-1} |x_{1,l} - x_{2,l}|
  \leq \rho_{n,l} \leq C |x_{1,l} - x_{2,l}|, \;
  \int_{B_{2 \rho_{n,l}}(y_{n,l})} e^{4 u_l} d V_g \leq \rho \quad n = 1, \dots,
  \ell_l.
\end{equation}
Reasoning as in the proof of Proposition \ref{p:cc} one finds
\begin{eqnarray*}
    \int_{B_{\rho_{n,l}}(y_{n,l})} e^{4 u_l} d V_g & \leq & C
   \int_{B_{\rho_{n,l}}(y_{n,l})} d V_g(x) \exp \left[ 4 \int_{M \setminus
   B_{2 \rho_{n,l}}(y_{n,l})} G(x,y) 2 k_l e^{4 u_l(y)} d V_g(y)
   \right] \\ & \times & \int_{B_{2 \rho_{n,l}}(y_{n,l})} \left(
   \frac{1}{|x - y|} \right)^{\frac{k_l \rho}{\pi^2}} e^{4 \ov{u}_l}
   d V_g(y).
\end{eqnarray*}
From \eqref{eq:G} and \eqref{eq:eded}, after some computation we
get
\begin{eqnarray}\label{eq:l1s3}
  \int_{B_{\rho_{n,l}}(y_{n,l})} e^{4 u_l} d V_g & \leq & C
  \int_{B_{\rho_{n,l}}(y_{n,l})}
  \left[ 4 \int_{M \setminus \left( B_{2 \rho_{n,l}}(y_{n,l}) \cup
  B_{\frac{C^{-1} |x_{1,l}-x_{2,l}|}{2}}(x_{1,l})\right)} G(x,y) 2
  k_l e^{4 u_l(y)} d V_g(y)
  \right] \nonumber \\ & \times & |x_{1,l} - x_{2,l}|^{-8+o_l(1)}
  |x_{1,l} - x_{2,l}|^{4-\frac{k_l \rho}{\pi^2}} e^{4 \ov{u}_l}
  d V_g(x) \\ & \leq & C \sup_{x \in B_{\rho_{n,l}}(y_{n,l})}
  \left[ 8 \int_{M \setminus \left(B_{2 \rho_{n,l}}(y_{n,l}) \cup
  B_{\frac{C^{-1} |x_{1,l}-x_{2,l}|}{2}}(x_{1,l})\right)}
  G(x,y) k_l e^{4 u_l(y)} d V_g(y) \right] \nonumber \\ & \times &
  |x_{1,l} - x_{2,l}|^{-\frac{k_l \rho}{\pi^2}+o_l(1)} e^{4 \ov{u}_l}.
  \nonumber
\end{eqnarray}
since $\rho_{n,l}$ is bounded from above by $C |x_{1,l} -
x_{2,l}|$.

On the other hand, if $\tilde{\e}$ and $R_{\tilde{\e}}$ are as in
Step 2, we also have
\begin{eqnarray*}
    u_l(x) & \geq & - C + \ov{u}_l + \int_{M \setminus B_{\frac{C^{-1}
  |x_{1,l}-x_{2,l}|}{2}}(x_{1,l})} G(x,y) 2 k_l e^{4 u_l(y)} d
  V_g(y) \\ & + & \int_{B_{R_{\tilde{\e}} r_{1,l}}(x_{1,l})} G(x,y)
  2 k_l e^{4 u_l(y)} d V_g(y), \qquad \hbox{ a.e. } x
  \in B_{\frac{C^{-1} |x_{1,l}-x_{2,l}|}{4}}(x_{1,l}) \setminus
  B_{2 R_{\tilde{\e}} r_{1,l}}(x_{1,l}).
\end{eqnarray*}
Reasoning as for \eqref{eq:eqeq3}, we then deduce that
\begin{eqnarray*}
    1 & \geq & \int_{B_{\frac{C^{-1} |x_{1,l}-x_{2,l}|}{4}}(x_{1,l})
    \setminus B_{2 R_{\tilde{\e}} r_{1,l}}(x_{1,l})} e^{4 u_l} d
    V_g \geq C^{-1} e^{4 \ov{u}_l} (R_{\tilde{\e}} r_{1,l})^{\tilde{\e} -
    4} \\ & \times & \inf_{z \in B_{\frac{C^{-1}
    |x_{1,l}-x_{2,l}|}{4}}(x_{1,l})} \left[ 8 \int_{M \setminus
    \left(B_{2 \rho_{n,l}}(y_{n,l}) \cup B_{\frac{C^{-1}
    |x_{1,l}-x_{2,l}|}{2}}(x_{1,l})\right)} G(z,y) k_l e^{4 u_l(y)}
    d V_g(y) \right].
\end{eqnarray*}
Now we notice that by \eqref{eq:rhoil} and \eqref{eq:G} one has
\begin{eqnarray*}
    |G(z,y) - G(x,y)| \leq C; \qquad x \in
    B_{\rho_{n,l}}(y_{n,l}), y \in M \setminus
    \left(B_{2 \rho_{n,l}}(y_{n,l}) \cup B_{\frac{C^{-1}
    |x_{1,l}-x_{2,l}|}{2}}(x_{1,l})\right), \\ \hbox{ and for }
    z \in B_{\frac{C^{-1} |x_{1,l}-x_{2,l}|}{4}}(x_{1,l}).
\end{eqnarray*}
From \eqref{eq:l1s3} and the last two formulas it follows that
$$
  \int_{B_{\rho_{n,l}}(y_{n,l})} e^{4 u_l} d V_g \leq C
  |x_{1,l} - x_{2,l}|^{-\frac{k_l \rho}{\pi^2}+o_l(1)} (R_{\tilde{\e}}
  r_{1,l})^{\tilde{\e} - 4} \to 0, \qquad \hbox{ as } l \to +
  \infty,
$$
since $\frac{r_{1,l}}{|x_{1,l} - x_{2,l}|} \to 0$ by
\eqref{eq:rhrh2}. Then the conclusion follows from \eqref{eq:eded}
and the fact that $B_{C |x_{1,l} - x_{2,l}|}(x_{1,l}) \setminus
\cup_{i = 1, \dots, h} B_{C^{-1} |x_{1,l} - x_{2,l}|}(x_{i,l})$ is
covered by a finite (and uniformly bounded) number of balls
$B_{\rho_{n,l}}(y_{n,l})$.
\end{pf}

\

\noindent Now we let
$$
  d_{1,l} = \inf \left\{ |x_{1,l} - x_{i,l}| \ \; : \; x_{i,l}
  \not\in X_{1,l} \right\}.
$$
Note that, by our definition of $X_{1,l}$, we have
$\frac{d_{1,l}}{|x_{1,l} - x_{2,l}|} \to + \infty$ as $l \to +
\infty$. We prove next the following result, which improves the
estimate in formula \eqref{eq:star2} to a larger set.

\begin{lem}\label{l:refin}
There exists $C > 0$ such that for $l$ large
\begin{equation}\label{eq:refin}
    \int_{B_{C^{-1} d_{1,l}}(x_{1,l})} e^{4 u_l} dV_g =
    \frac{8 \pi^2}{k_0} card(X_{1,l}) + o_l(1).
\end{equation}
\end{lem}

\begin{pf}
The proof follow closely the arguments of Proposition \ref{p:sbu},
hence we will be sketchy. We use the same notation as in Section
\ref{s:sbu} for the functions $(W_l)_l$ and the annuli $A_{r,l}$,
except for the fact that now we take $x_{1,l}$ as centers, hence
replacing the points $x_l$.

First of all we notice that, by the arbitrarity of $C$ in Lemma
\ref{l:step31}, there exists $Z_l \to + \infty$ such that
\begin{equation}\label{eq:st3ann}
    \int_{B_{e^{4 Z_l} |x_{1,l} - x_{2,l}|}(x_{1,l}) \setminus
  B_{C |x_{1,l} - x_{2,l}|}(x_{1,l})} e^{4 u_l} d V_g
  \to 0 \qquad \hbox{ as } l \to + \infty.
\end{equation}
Using the Jensen inequality in the annulus $B_{e^{4 Z_l} |x_{1,l}
- x_{2,l}|}(x_{1,l}) \setminus B_{e^{Z_l} |x_{1,l} -
x_{2,l}|}(x_{1,l})$, it follows that
\begin{equation}\label{eq:unif2}
    \sup_{z \in [ Z_l + \log |x_{1,l} - x_{2,l}|, 4 Z_l + \log
    |x_{1,l} - x_{2,l}|} \left( z + W_l(z) \right) \to - \infty
    \qquad \hbox{ as } l \to + \infty.
\end{equation}
Our next goal is to prove that also
\begin{equation}\label{eq:cllast}
  W'_l(z) = - 2 card(X_{1,l}) + o_l(1); \qquad \hbox{ for }
  z \in [ 2 Z_l + \log |x_{1,l} - x_{2,l}|, 3 Z_l + \log |x_{1,l} -
  x_{2,l}| ].
\end{equation}
In order to show this, we notice that by the second formula in
Remark \ref{r:bdW'l} and by some manipulation (reasoning as in the
proof of Lemma \ref{l:system}), there holds
$$
  W'_l(z) = \frac{\int_{A_{r,l}} \frac{\partial u_l}{\partial z}
  f e^{4z}}{\int_{A_{r,l}} f e^{4z}} + O(e^{2z}); \qquad \hbox{ for }
  z \in [ Z_l + \log |x_{1,l} - x_{2,l}|, 4 Z_l + \log |x_{1,l} -
  x_{2,l}| ], r = e^z.
$$
Using the Green's representation formula we obtain
\begin{eqnarray*}
  \frac{\partial u}{\partial r}(x) & = & \int_{B_{e^{Z_l} |x_{1,l} -
  x_{2,l}|}(x_{1,l})} \frac{\partial_x G(x,y)}{\partial r} f_l(y) d
  V_g(y) + \int_{M \setminus B_{e^{4 Z_l} |x_{1,l} -
  x_{2,l}|}(x_{1,l})} \frac{\partial_x G(x,y)}{\partial r} f_l(y) d
  V_g(y) \\ & + & \int_{B_{e^{4 Z_l} |x_{1,l} - x_{2,l}|}(x_{1,l})
  \setminus B_{e^{Z_l} |x_{1,l} - x_{2,l}|}(x_{1,l})}
  \frac{\partial_x G(x,y)}{\partial r} f_l(y) d V_g(y).
\end{eqnarray*}
From \eqref{eq:derG2}, Lemma \ref{l:step31} and \eqref{eq:st3ann}
it follows that, for $Z_l \to + \infty$ sufficiently slowly
$$
  \int_{B_{e^{Z_l} |x_{1,l} - x_{2,l}|}(x_{1,l})} \frac{\partial_x
  G(x,y)}{\partial r} f_l(y) d V_g(y) = - \frac{2
  card(X_{1,l})}{|x - x_{1,l}|} + o_l(1).
$$
Also, reasoning as in the proof of Lemma \ref{l:lodlp} and using
\eqref{eq:st3ann} one finds that
$$
  \left| \int_{A_{r,l}} dx \int_{B_{e^{4 Z_l} |x_{1,l} - x_{2,l}|}(x_{1,l})
  \setminus B_{e^{Z_l} |x_{1,l} - x_{2,l}|}(x_{1,l})}
  \frac{\partial_x G(x,y)}{\partial r} f_l(y) d V_g(y) \right| = o(1) |x -
  x_{1,l}|^3.
$$
Finally, since $Z_l \to + \infty$ one also finds that
$$
  \int_{M \setminus B_{e^{4 Z_l} |x_{1,l} - x_{2,l}|}(x_{1,l})}
  \frac{\partial_x G(x,y)}{\partial r} f_l(y) d V_g(y) =
  o_l(1) \frac{1}{|x - x_{1,l}|}.
$$
Recalling that $\frac{\partial u_l}{\partial z} = r \frac{\partial
u_l}{\partial r}$, with $r = dist(x, x_{1,l})$, the last three
formulas yield \eqref{eq:cllast}.

Now, for $\g \in (1, 2)$ we consider the following sequence of
functions
$$
  h^\g_{l}(z) = - \g (z - \log |x_{1,l} - x_{2,l}| - 2 Z_l) +
  W_l(\log |x_{1,l} - x_{2,l}| + 2 Z_l).
$$
Exactly as in the proof of Proposition \ref{p:sbu} one can show
that
$$
  W_l(z) \leq h_{l}^\g(z); \qquad \qquad z \in \left[
  \log |x_{1,l} - x_{2,l}| + 2 Z_l, \log d_{1,l} - C_\g \right].
$$
As above, we define
$$
  \hat{r}_{l,j} = \frac 45 e^{2 Z_l} |x_{1,l} - x_{2,l}| \left( \frac 75
  \right)^j ; \qquad \qquad \left( \frac 75 \right)^{j_l} \in \left(
  \frac 54 \frac{d_{1,l}}{C e^{2 Z_l} |x_{1,l} - x_{2,l}|}, \frac 52
  \frac{d_{1,l}}{C e^{2 Z_l} |x_{1,l} - x_{2,l}|} \right).
$$
and we obtain
\begin{eqnarray*}
    \int_{A'_{\hat{r}_{l,j},l}} e^{4 u_l} dV_g \leq C
    |A_{\hat{r}_{l,j},l}| e^{4 \ov{u}_{l,\hat{r}_{l,j}}} \leq C
    \hat{r}_{l,j}^4 e^{4 W_l(\log \hat{r}_{l,j})} \leq C
    \hat{r}_{l,j}^4 e^{4 h_{l}^\g(\log \hat{r}_{l,j})}; \qquad
    j = 1, \dots, j_l.
\end{eqnarray*}
From the expression of $h^\g_{l}$ and \eqref{eq:unif2} we deduce
\begin{eqnarray*}
    \hat{r}_{l,j}^4 e^{4 h_{l}^\g(\log \hat{r}_{l,j})} & \leq &
    C \hat{r}_{l,j}^4 \exp \left[ 4 \left( - \g \left( \log \hat{r}_{l,j}
    - \log |x_{1,l} - x_{2,l}| - 2 Z_l \right) + W_l(\log |x_{1,l} - x_{2,l}| + 2 Z_l) \right)
    \right] \\ & \leq & o_l(1) \hat{r}_{l,j}^4 \exp \left[ - 4 \g \log \hat{r}_{l,j}
    + 4 (\g - 1) \log |x_{1,l} - x_{2,l}| + 8 (\g-2) Z_l \right]
    \\ & \leq & o_l(1) \left( \frac{|x_{1,l} - x_{2,l}|}{\hat{r}_{l,j}}
    \right)^{4(\g-1)} = o_l(1) \left( \frac 57 \right)^{4(\g-1)j}.
\end{eqnarray*}
As before we then find
$$
  \int_{B_{\frac{d_{1,l}}{C}}(x_l) \setminus B_{C |x_{1,l} - x_{2,l}|}(x_l)}
  e^{4 u_l} dV_g \leq o_l(1) \sum_{j=0}^{\infty} \left( \frac 57
\right)^{4(\g-1)j} \to 0,
$$
This formula, joint with \eqref{eq:star2}, yields the conclusion
of the Lemma.
\end{pf}

\

\noindent The proof of Step 3 follows from the arguments of Lemmas
\ref{l:step31}, \ref{l:refin}, repeating the procedure for all the
clusters of the points of $\{ x_{1,l}, \dots, x_{j,l}\} \setminus
X_{1,l}$.

\

\noindent The proof of the theorem is now an easy consequence of
\eqref{eq:noul} and \eqref{eq:star4}, since $k_0$ is not an
integer multiple of $8 \pi^2$.

\subsection{The case $k_0 < 8 \pi^2$}\label{ss:pnot>0}

In this final subsection we consider the cases in which $P_g$
possesses some negative eigenvalues and $k_0 < 8 \pi^2$.  We prove
first the following result, which regards boundedness of the
$V$-component of sequences of solutions.

\begin{lem}\label{l:bdV}
Suppose $P_g$ possesses some negative eigenvalues, and suppose
that $ker P_g = \{constants\}$. Let $(u_l)_l \subseteq H^2(M)$ be
a sequence satisfying \eqref{eq:pl}-\eqref{eq:noul}. Let us write
$u_l = \hat{u}_l + \tilde{u}_l$ with $\hat{u}_l \in V$ and
$\tilde{u}_l \perp V$, where $V$ denotes the direct sum of the
negative eigenspaces of $P_g$. Then there holds
$$
  \|\hat{u}_l\|_{H^2(M)} \leq C,
$$
for some positive constant $C$ independent of $l$.
\end{lem}

\begin{pf}
Let $\hat{v}_1, \dots, \hat{v}_{\ov{k}}$ be as in
\eqref{eq:hatv1k}. Then, by standard elliptic regularity theory,
each $\hat{v}_i$ is smooth on $M$. Testing \eqref{eq:pl} on
$\hat{u}_l$ we obtain
$$
  \langle P_g \hat{u}_l, \hat{u}_l \rangle + 4 \int_M Q_l
  \hat{u}_l dV_g + 4 k_l \int_M e^{4 u_l} \hat{u}_l dV_g = 0.
$$
Using \eqref{eq:noul}, the fact that on $V$ the $L^\infty$-norm is
equivalent to the $H^2$-norm, and the Poincar\'e inequality, from
the last formula we deduce that
$$
   - \langle P_g \hat{u}_l, \hat{u}_l \rangle = O(1) \|\hat{u}_l\|_{H^2(M)}.
$$
Since $P_g$ is negative-definite on $V$, the conclusion follows.
\end{pf}

\

\noindent Next, we consider separately the following three
possibilities, one of which will always occur for $k_0 < 8 \pi^2$
and for $l$ sufficiently large.

\

\noindent {\bf Case 1: $k_l < 0$.} First of all, using the Jensen
inequality we find immediately that $\ov{u}_l \leq C$, for some
constant $C$ independent of $l$. Then, multiplying \eqref{eq:pl}
by $u_l$ and integrating on $M$, using the Poincar\'e inequality
and Lemma \ref{l:bdV}, we find
\begin{eqnarray*}
  \langle P_g u_l, u_l \rangle & = & 2 k_l \int_M e^{4 u_l} u_l d
  V_g - 2 k_l \ov{u}_l + O \left( \langle P_g u_l, u_l \rangle^{\frac 12}
  \right) + C \\ & \leq & C + (- 2 k_l) \ov{u}_l + O \left( \langle P_g
  u_l, u_l \rangle^{\frac 12} \right) \leq C + O \left( \langle P_g u_l,
  u_l \rangle^{\frac 12} \right).
\end{eqnarray*}
Again by Lemma \ref{l:bdV}, this implies uniform bounds on $\|u_l
- \ov{u}_l\|$ and hence, by \eqref{eq:adaneg}, uniform $L^p$
bounds on $e^{4 u_l}$ for any $p > 1$. Then the conclusion follows
from standard elliptic regularity results.

\

\noindent {\bf Case 2: $0 \leq k_l 2 \leq \pi^2$.} Since we are
assuming \eqref{eq:noul}, we easily see that the alternative
\eqref{eq:conc} in Proposition \ref{p:cc} cannot occur. Therefore,
reasoning as in the previous case, we obtain uniform $L^p$ bounds
on $e^{4 u_l}$ for some $p > 1$.

\

\noindent {\bf Case 3: $2 \pi^2 \leq k_l < \frac{1}{2} (k_0 +
8\pi^2) < 8 \pi^2$.} In this case it is $k_0 > 0$. Assuming
$(u_l)_l$ unbounded, Proposition \ref{p:hwlb} applies, and
\eqref{eq:intalrlu} gives a contradiction to \eqref{eq:noul},
since $k_0 < 8 \pi^2$.

\end{document}